\documentclass[12pt,a4paper]{article}

\usepackage[maxnames=4,minnames=3,maxbibnames=99]{biblatex}
\bibliography{bib}
\usepackage{amsfonts}
\usepackage{amsmath}
\usepackage{amssymb}
\usepackage{amsthm}
\usepackage{mathrsfs}
\usepackage{yfonts}
\usepackage{pifont}
\usepackage{marvosym}
\usepackage[dvips]{graphicx}
\usepackage[all]{xy}
\usepackage{stmaryrd}
\usepackage{color}
\usepackage{tikz}
\usepackage[margin=0.8in]{geometry}
\usetikzlibrary{arrows,matrix}

\newcommand{\RR}{\mathbb{R}}

\newcommand{\XX}{\mathbb{X}}
\newcommand{\YY}{\mathbb{Y}}
\newcommand{\VV}{\mathbb{V}}

\def\colim{\mathop{\rm colim}\nolimits}

\theoremstyle{definition}
 \theoremstyle{plain}

\newtheorem{Theorem}{Theorem}[section]
\newtheorem{Lemma}[Theorem]{Lemma}

\newtheorem{Remark}[Theorem]{Remark}

\newtheorem*{Lemma*}{Lemma}

\theoremstyle{definition}

\newtheorem{Definition}[Theorem]{Definition}
\newtheorem*{Definition*}{Definition}
\newtheorem{Example}[Theorem]{Example}
\newtheorem*{Example*}{Example}

\usepackage{graphicx}

 \usepackage{graphicx}

 \usepackage{graphicx}

 \usepackage{graphicx}

\title{Alexander Duality for Parametrized Homology
}

\author{
Sara Kali\v{s}nik\vspace{1ex}\\
Melvin and Joan Lane Stanford Graduate Fellow\\
Department of Mathematics\\
Stanford University\\
\texttt{kalisnik@stanford.edu}\vspace{1ex}\\
}
\date{\vspace{-5ex}}
\begin{document}
\maketitle
\begin{abstract}
This paper extends Alexander duality to the setting of parametrized homology. Let  $X \subset \RR^n \times \RR$ with $n\geq 2$ be a compact set satisfying certain conditions, let $Y = (\RR^n \times \RR) \setminus X$, and let $p$ be the projection onto the second factor. Both $X$ and $Y$ are parametrized spaces with respect to the projection.  The parametrized homology is a variant of zigzag persistent homology that measures how the homology of  the level sets of the space changes as we vary the parameter. We show that if $(X, p|_X)$ has a well-defined parametrized homology, then  the pair $(Y, p|_Y)$ has a well-defined reduced parametrized homology. We also establish a relationship between the parametrized homology of $(X, p|_X)$ and the reduced parametrized homology of $(Y, p|_Y)$.
\end{abstract}


\section{Introduction}

It is well-accepted that topological techniques can be useful for understanding high dimensional data. Computational topologists view data as finite metric spaces, build different complexes on the points (\v{C}ech, Vietoris-Rips), and analyze the topology of those objects to infer the topology of the data. This process is motivated by the nerve theorem in algebraic topology, which claims that given a covering of the space with balls, the \v{C}ech complex associated with this covering is homotopy equivalent to the space.

Building any of these different complexes requires a choice of parameter,  such as the radius of the balls in the case of the \v{C}ech complex. The idea of persistent homology is to let the parameter value vary while tracking the births and deaths of topological features. The output  is a persistence diagram that measures the significance of a topological feature.

To ensure that the homology would change at finitely many values, scholars have imposed various restrictions (such as assuming a function to be Morse or a space to be compact). Chazal et al.\,\cite{Thestructureandstability} avoid these restrictions and define persistence diagrams in a wider variety of situations. This approach can also be used to define the levelset zigzag persistence \cite{Zigzagpersistenthomologyandreal}  more broadly with the outcome of the \emph{parametrized homology} \cite{levelsetmeasure}.

We are not only interested in determining relationships between homology groups, but also between persistence diagrams. For example, Cohen-Steiner et al.\,\cite{Extendingpersistence} deal with the question of extending the Poincar\'e duality to persistence diagrams.
Another classical theorem in algebraic topology is Alexander duality, which asserts a relationship between the homology groups of a locally contractible compact space and its complement. 

The goal of this paper is to extend this theorem to the setting of parametrized homology.  Let  $X \subset \RR^n \times \RR$ with $n\geq 2$ be a compact set, let $Y = (\RR^n \times \RR) \setminus X$,
and let $p$ be the projection onto the second factor. We assume that level sets $p^{-1}(a) \cap X$ for $a\in \RR$, and slices $p^{-1}([a,b]) \cap X$ for $a<b$ are locally contractible. We show that if $(X, p|_X)$ has a well-defined parametrized homology, then  the pair $(Y, p|_Y)$ has a well-defined reduced parametrized homology. More specifically, we show that the reduced parametrized homology of $(Y, p|_Y)$ in dimension $n-j-1$ is equal to the parametrized homology of  $(X, p|_X)$ in dimension $j$ for $j =0, \ldots, n-1$. While it may not be immediately obvious, the complement is well-behaved as a consequence of Alexander duality. We also establish a duality in terms of how homological features perish. If a $j$-dimensional homology cycle in $X$ is killed at the parameter value $p$, then there is a corresponding $(n-j-1)$-dimensional homology cycle in $Y$ that ceases to exist beyond $p$ and vice versa (see Section 2.2). Our theorem includes cases $(X, p|_X)$ where:
\begin{itemize}
\item
 $X$ is a compact submanifold of $\RR^n\times \RR$ (with or without boundary) and $p|_X$ is Morse; or, more generally, when $(X, p|_X)$ is of Morse-type \cite{Zigzagpersistenthomologyandreal} and $X$ is compact with all the slices and level sets locally contractible; 
 \item
$X$ is a finite simplicial complex and $p|_X$ is a piecewise-linear map;
\item
$X$ is a semialgebraic subset of $\RR^n\times \RR$.
 \end{itemize}

Edelsbrunner and Kerber have expanded Alexander duality to extended persistence diagrams \cite{AlexanderDuality}. However, there is a difference in our approaches to  Alexander duality. Edelsbrunner and Kerber consider Alexander duality to be a statement about two complementary subsets of the sphere that intersect in a $n$-manifold. By contrast, we consider it to be a statement about a compact
subset of a Euclidean space (or a sphere) and its complement. An advantage of our approach is that it allows us to generalize
Alexander duality directly to an appropriate parametrized version,
starting with a compact parametrized space that satisfies certain conditions. 

One possible application of this research is in sensor networks \cite{sensor}. The classic Alexander duality can be used to find gaps in static networks. However, since we are interested in time-varying networks, we develop a parametrized version of this theorem. While we can observe the space covered by the sensors within a time-varying network, this method gives us knowledge of the uncovered regions.

\section{Background}
Throughout this paper we work with homology and cohomology with coefficients in a field $\mathbf k$. So  
$\operatorname{\rm H}_j(X)$ always means $\operatorname{\rm H}_j(X, \mathbf k)$ and $\operatorname{\rm H}^j(X)$ means $\operatorname{\rm H}^j(X, \mathbf k)$. 

\subsection{Persistence diagrams}

Persistence is commonly described either as a multiset of intervals (a barcode) or as a multiset of points in the half plane (a persistence diagram). However, these descriptions do not distinguish between different types of intervals (for example, $[p, q]$, $[p, q)$, $(p, q]$ and $(p, q)$). In order to capture the homology of spaces parametrized  over $\RR$, Chazal et al.\,\cite{Thestructureandstability} introduced decorated real numbers. We can represent every decorated point as a point in the half plane with a tick specifying the decoration:
We adopt the following notation: 
$$
\begin{array}{clclllll}
[ p, q) & \textrm{is written} &(p^-, q^-) & \textrm{and drawn} & \scalebox{0.13}{\definecolor{c090000}{RGB}{9,0,0}
\definecolor{cff8080}{RGB}{255,128,128}

\begin{tikzpicture}[y=0.80pt, x=0.8pt,yscale=-1, inner sep=0pt, outer sep=0pt]
\path[cm={{-1.0,0.0,0.0,-1.0,(-51.15,-104.83783)}},draw=c090000,fill=cff8080,line
  cap=rect,miter limit=4.00,line width=4.560pt]
  (-130.0000,-93.3500)arc(0.000:180.000:30.000)arc(-180.000:0.000:30.000) --
  cycle;

\path[draw=black,line join=miter,line cap=butt,miter limit=4.00,line
  width=4.560pt] (88.8500,8.5122) -- (58.8500,38.5122);

\end{tikzpicture}}\\
\left[ p, q \right] & \textrm{is written} & (p^-, q^+) & \textrm{and drawn}& \scalebox{0.13}{\definecolor{c090000}{RGB}{9,0,0}
\definecolor{cff8080}{RGB}{255,128,128}

\begin{tikzpicture}[y=0.80pt, x=0.8pt,yscale=-1, inner sep=0pt, outer sep=0pt]
\path[cm={{-1.0,0.0,0.0,1.0,(-51.15,101.86218)}},draw=c090000,fill=cff8080,line
  cap=rect,miter limit=4.00,line width=4.560pt]
  (-130.0000,-93.3500)arc(0.000:180.000:30.000)arc(-180.000:0.000:30.000) --
  cycle;

\path[draw=black,line join=miter,line cap=butt,miter limit=4.00,line
  width=4.560pt] (88.8500,-11.4878) -- (58.8500,-41.4878);

\end{tikzpicture}}\\
(p, q) &  \textrm{is written} & (p^+, q^-)  & \textrm{and drawn} & \scalebox{0.13}{\definecolor{c090000}{RGB}{9,0,0}
\definecolor{cff8080}{RGB}{255,128,128}

\begin{tikzpicture}[y=0.80pt, x=0.8pt,yscale=-1, inner sep=0pt, outer sep=0pt]
\path[cm={{1.0,0.0,0.0,-1.0,(248.85,-104.83783)}},draw=c090000,fill=cff8080,line
  cap=rect,miter limit=4.00,line width=4.560pt] (-130.0000,-93.3500) .. controls
  (-130.0000,-76.7815) and (-143.4315,-63.3500) .. (-160.0000,-63.3500) ..
  controls (-176.5685,-63.3500) and (-190.0000,-76.7815) .. (-190.0000,-93.3500)
  .. controls (-190.0000,-109.9185) and (-176.5685,-123.3500) ..
  (-160.0000,-123.3500) .. controls (-143.4315,-123.3500) and
  (-130.0000,-109.9185) .. (-130.0000,-93.3500) -- cycle;

\path[draw=black,line join=miter,line cap=butt,miter limit=4.00,line
  width=4.560pt] (108.8500,8.5122) -- (138.8500,38.5122);

\end{tikzpicture}}\\
(p, q ] &  \textrm{is written}& (p^+, q^+) & \textrm{and drawn} &   \scalebox{0.13}{\definecolor{c090000}{RGB}{9,0,0}
\definecolor{cff8080}{RGB}{255,128,128}

\begin{tikzpicture}[y=0.80pt, x=0.8pt,yscale=-1, inner sep=0pt, outer sep=0pt]
\path[shift={(248.85,101.86218)},draw=c090000,fill=cff8080,line cap=rect,miter
  limit=4.00,line width=4.560pt] (-130.0000,-93.3500) .. controls
  (-130.0000,-76.7815) and (-143.4315,-63.3500) .. (-160.0000,-63.3500) ..
  controls (-176.5685,-63.3500) and (-190.0000,-76.7815) .. (-190.0000,-93.3500)
  .. controls (-190.0000,-109.9185) and (-176.5685,-123.3500) ..
  (-160.0000,-123.3500) .. controls (-143.4315,-123.3500) and
  (-130.0000,-109.9185) .. (-130.0000,-93.3500) -- cycle;

\path[draw=black,line join=miter,line cap=butt,miter limit=4.00,line
  width=4.560pt] (108.8500,-11.4878) -- (138.8500,-41.4878);

\end{tikzpicture}}\\
\end{array}
$$
We require $p<q$, except for the one point interval $[p, p]$. The notation for an arbitrary interval is $(p^*, q^*)$.

Let  $R = [a, b ] \times [c, d]$, where $a < b \leq c < d$, be a rectangle. Let $(p^*, q^*)$ be a decorated point. Then $(p^*, q^*) \in R$ if $[b, c] \subset (p^*, q^*) \subset (a, d)$. This happens exactly when the point $(p, q)$ and its decoration tick are contained in the closed rectangle $R$:\begin{center}
\includegraphics[scale=0.45]{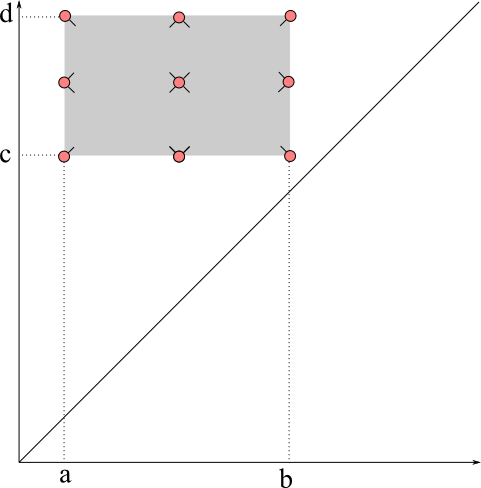}
\end{center}
A locally finite multiset of decorated points in the half plane is called a decorated persistence diagram. 

Chazal et al.\,\cite{Thestructureandstability} introduce a new approach for expressing persistence that is especially well-suited for a continuous parameter. The intuition is that if we know how many points of the diagram are contained in each rectangle in the half plane, then we know the diagram itself. Counting the points in the rectangles leads to the introduction of r-measures. 

We work in the open half plane $\mathscr{H} =\{(p, q) \in \RR^2 \, \mid \, p<q \}$ throughout this paper, so we state the relevant results for this case only. 
\begin{Definition} 
The set of rectangles in $\mathscr{H}$ is
$$
\operatorname{\rm Rect}(\mathscr{H}) = \{[a, b]\times [c,d]\subset \mathscr{H} \, |\, a<b<c<d\}.
$$
A rectangle measure or r-measure on $\mathscr{H}$ is a function 
$$
\mu \colon \operatorname{\rm Rect}(\mathscr{H})  \to \{0, 1, 2, 3, \ldots \}\cup \{\infty \}
$$
that is additive under vertical and horizontal splitting, meaning that $\mu(R) = \mu(R_1) + \mu(R_2)$ whenever
\begin{center}
\fbox{$R$} = \fbox{$R_1$}\fbox{$R_2$} or \fbox{$R$} = \begin{tabular}{|c|}
\hline
$R_1$ \\
\hline
$R_2$\\
\hline
\end{tabular}  .
\end{center} 

\end{Definition}

Finite r-measures correspond exactly to decorated persistence diagrams \cite{Thestructureandstability}.
\begin{Theorem}[The Equivalence Theorem]\label{equiv}
There is a bijective correspondence between:
\begin{itemize}
\item
Finite r-measures $\mu$ on $\mathscr{H} $. Here `finite' means that $\mu(R)<\infty$ for every $R\in \operatorname{\rm Rect}(\mathscr{H}) $.
\item
Locally finite multisets $A$ of decorated points in $\mathscr{H}$.
Here `locally finite' means that $\operatorname{\rm card}(A|_R)<\infty$ for every $R\in \operatorname{\rm Rect}(\mathscr{H})$.
\end{itemize}
The measure $\mu$ corresponding to a multiset $A$ satisfies the formula
$$
\mu(R) = \operatorname{\rm card}(A|_R)
$$
for every $R\in \operatorname{\rm Rect}(\mathscr{H})$.
\end{Theorem}

From the proof of the equivalence theorem, we get the locally finite multiset $A$ determined by the measure $\mu$ by computing multiplicities. The multiplicity of $(p^*, q^*)$ with respect to $\mu$ is
$$
\textrm{m}_\mu(p^*, q^*) = \textrm{min}\{\mu(R)\, \mid\, (p^*, q^*) \in R, R \in \operatorname{\rm Rect}(\mathscr{H})\}.
$$
Alternatively, we can pick a nested sequence $R_1 \supset R_2 \supset R_3 \supset \ldots$  of closed rectangles that contain $(p^*, q^*)$ such that $\cap_n R_n = (p, q)$ (see Figure 1). Then
$$
\textrm{m}_\mu(p^*, q^*) = \lim_{n\to \infty} \mu(R_n).
$$
\begin{figure}[h!]\label{mult}
\begin{center}
\includegraphics[scale=0.1]{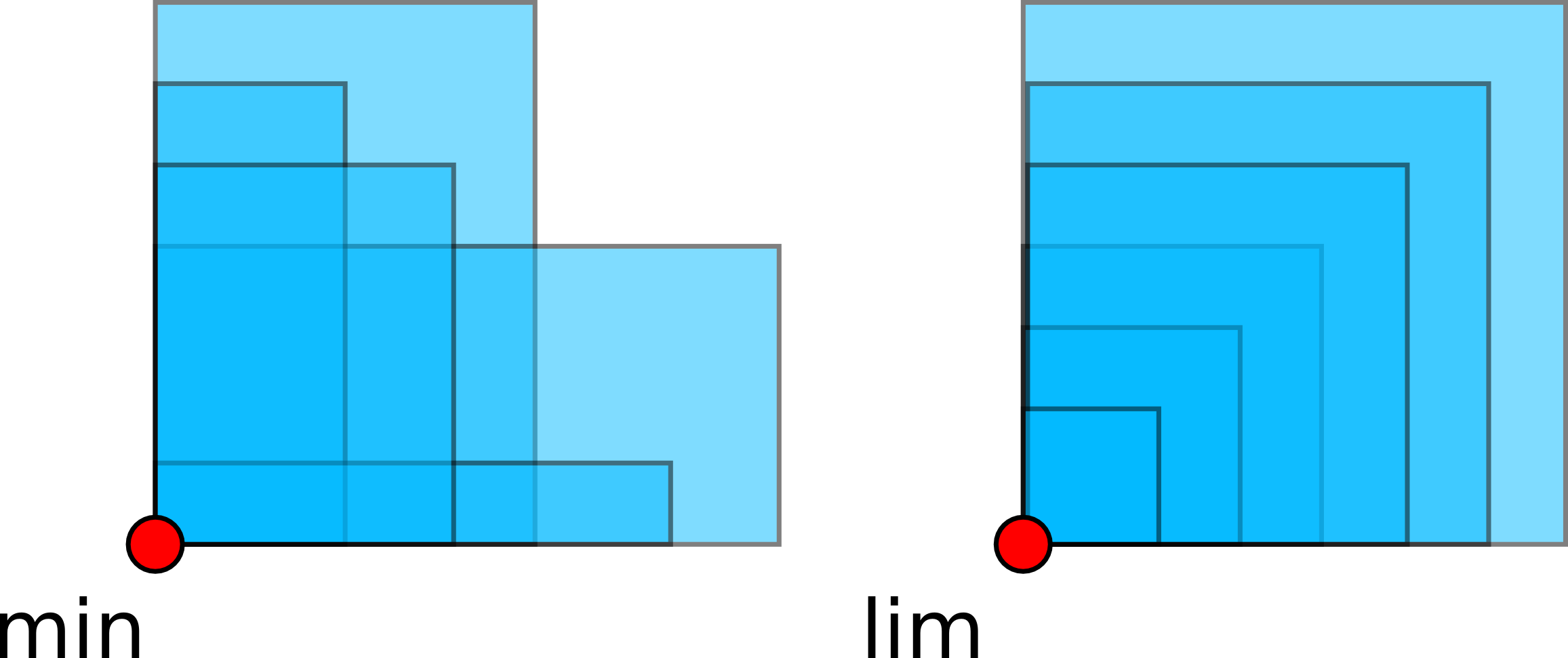}
\caption{Computing multiplicities.}
\end{center}
\end{figure}
We use these formulas to determine the decorated diagram in Example~\ref{example}.

\subsection{Parametrized homology}
Carlsson et al.\,introduced levelset zigzag persistence \cite{Zigzagpersistenthomologyandreal} and presented it in the measure context as parametrized homology \cite{levelsetmeasure}. A parametrized space is a pair $\XX = (X, p )$, where $X$ is a topological space, and $p \colon X \to \RR$ is a continuous function. We are interested in how the homology changes as the parameter varies. The function $p$ defines levelsets $X_a = p^{-1}(a)$ and slices $X_a^b = p^{-1}([a,b])$.

Given a rectangle $R = [a,b]\times [c,d]$ with $-\infty <a<b<c<d < \infty$, the aim is to count the homological features of $\XX$ that persist over the closed interval $[b, c]$, but not over the open interval $(a, d)$. Consider the following diagram of spaces and inclusion maps.
\begin{center}
\begin{tikzpicture}[xscale=1,yscale=1, font=\small]
\draw (1,1) node(11) {$X_a^b$} (3,1) node(31) {$X_b^c$} (5,1) node(51) {$X_c^d$} ;
\draw (0,0) node(00){$X_a$} (2,0) node(20){$X_b$} (4,0) node(40){$X_c$} (6,0) node(60){$X_d$}  ;

\draw[->] (00) -- (11); 
\draw[->] (20) -- (11);
\draw[->] (20) -- (31);
\draw[->] (40) -- (31);
\draw[->] (40) -- (51);
\draw[->] (60) -- (51);  

\end{tikzpicture}
\end{center}
Apply the $j$-dimensional homology functor $\operatorname{\rm H}_j$ to obtain:
\begin{center}
\begin{tikzpicture}[xscale=1,yscale=1, font=\small]
\draw (1,1) node(11) {$\operatorname{\rm H}_j(X_a^b)$} (3,1) node(31) {$\operatorname{\rm H}_j(X_b^c)$} (5,1) node(51) {$\operatorname{\rm H}_j(X_c^d)$} ;
\draw (0,0) node(00){$\operatorname{\rm H}_j(X_a)$} (2,0) node(20){$\operatorname{\rm H}_j(X_b)$} (4,0) node(40){$\operatorname{\rm H}_j(X_c)$} (6,0) node(60){$\operatorname{\rm H}_j(X_d)$.}  ;

\draw[->] (00) -- (11); 
\draw[->] (20) -- (11);
\draw[->] (20) -- (31);
\draw[->] (40) -- (31);
\draw[->] (40) -- (51);
\draw[->] (60) -- (51);  

\end{tikzpicture}
\end{center}
Such a diagram of vector spaces and maps between them is called a zigzag module \cite{Zigzagpersistence}. It can be viewed as a representation of a quiver of type $A_7$.
We denote this quiver representation by $\operatorname{\rm H}_j(\XX_{\{a, b, c, d\}})$. It is decomposable by Gabriel's theorem \cite{quiver}. There are four types of indecomposable summands that meet $b$ and $c$, but not $a$ and $d$. By counting each of these summands, we get four quantities presented in the notation introduced by Chazal et al.\,\cite{Thestructureandstability}:
$$\begin{array}{rclcccccc}
{}_j\mu_\XX^{{ \backslash\! \backslash}}(R) &= &\langle \scalebox{0.065}{\input{dd.tex}}\, |\, \operatorname{\rm H}_j(\XX_{\{a, b, c, d\}}) \rangle \\
{}_j\mu_\XX^{\vee}(R)& = &\langle \scalebox{0.065}{\input{du.tex}}\, |\, \operatorname{\rm H}_j(\XX_{\{a, b, c, d\}}) \rangle\\
{}_j\mu_\XX^{\wedge}(R) &= &\langle \scalebox{0.065}{\input{ud.tex}}\, |\, \operatorname{\rm H}_j(\XX_{\{a, b, c, d\}}) \rangle\\
{}_j\mu_\XX^{^{/\!/}}(R) &= &\langle \scalebox{0.065}{\input{uu.tex}}\, |\, \operatorname{\rm H}_j(\XX_{\{a, b, c, d\}}) \rangle.
\end{array}
$$
Here $\langle \scalebox{0.065}{\input{dd.tex}}\, |\, \operatorname{\rm H}_j(\XX_{\{a, b, c, d\}}) \rangle$ denotes the number of times the summand $\scalebox{0.065}{\input{dd.tex}}$  appears in the interval decomposition of $\operatorname{\rm H}_j(\XX_{\{a, b, c, d\}})$. For the sake of simplicity we write \scalebox{0.065}{\input{dd.tex}} instead of $0 \to \mathbf k \leftarrow \mathbf k \to \mathbf k \leftarrow \mathbf k \to 0 \leftarrow 0$, where the maps $\mathbf k \to \mathbf k$ are identities and the other maps are 0.

Suppose these four quantities are finite r-measures. By the equivalence theorem each determines a decorated persistence diagram.  Let $\textrm{Dgm}_j^{*}(\XX)$  be the diagram determined by ${}_j\mu^*$. 

\begin{figure}[h!]
\begin{center}
 \scalebox{0.6}{\input{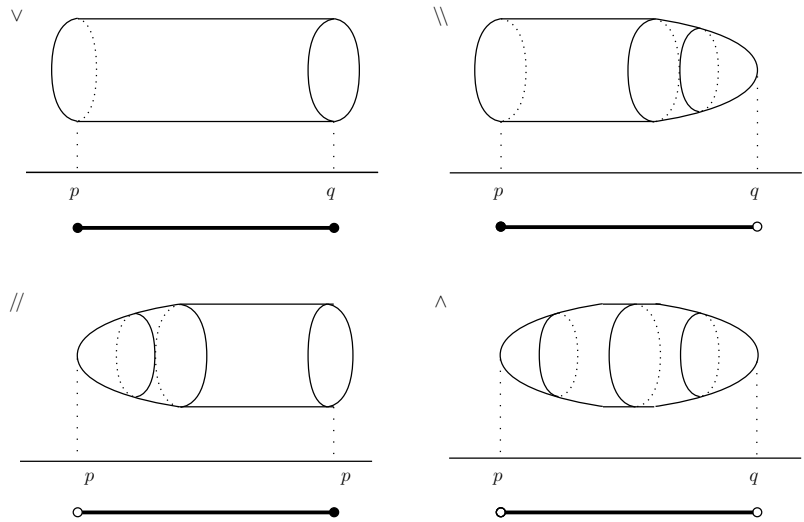}}
\caption{\label{perish} The 1-dimensional cycle on the upper left ceases to exist beyond both endpoints, whereas that on the upper right ceases to exist beyond $p$ and and is spanned by a disc at $q$.}
\end{center}
\end{figure}

These four diagrams demonstrate how homological features perish (whether $j$-dimensional cycles are \emph{killed} in 
homology by $(j+1)$-dimensional chains or whether they \emph{cease to exist}):
\begin{itemize}
\item
$\textrm{Dgm}_j^{\backslash \! \backslash}(\XX)$ contains decorated points $(p^*, q^*)$ corresponding to homology $j$-cycles that cease to exist beyond $p$, and are killed at $q$;
\item
$\textrm{Dgm}_j^{\vee}(\XX)$ contains decorated points $(p^*, q^*)$ corresponding to homology  $j$-cycles that cease to exist beyond both endpoints;
\item
$\textrm{Dgm}_j^{\wedge}(\XX)$ contains decorated points $(p^*, q^*)$ corresponding to homology  $j$-cycles that are killed at both endpoints;
\item
$\textrm{Dgm}_j^{/ \! /}(\XX)$ contains decorated points $(p^*, q^*)$ corresponding to homology $j$-cycles that are killed $p$ and cease to exist beyond $q$.
\end{itemize}
Figure~\ref{perish} shows examples of each type of homological feature discussed above.

\begin{Remark}
When $X$ is a compact manifold and $p$ is Morse, the four decorations correspondend exactly with how features perish at endpoints  \cite{levelsetmeasure}:
\begin{center}
\begin{tabular}{ccccccc}
${}_*\mu_\XX^{{ \backslash\! \backslash}} $  &
 \scalebox{0.2}{\definecolor{c090000}{RGB}{9,0,0}
\definecolor{cff8080}{RGB}{255,128,128}

\begin{tikzpicture}[y=0.80pt, x=0.8pt,yscale=-1, inner sep=0pt, outer sep=0pt]
\path[cm={{-1.0,0.0,0.0,-1.0,(-51.15,-104.83783)}},draw=c090000,fill=cff8080,line
  cap=rect,miter limit=4.00,line width=4.560pt]
  (-130.0000,-93.3500)arc(0.000:180.000:30.000)arc(-180.000:0.000:30.000) --
  cycle;

\path[draw=black,line join=miter,line cap=butt,miter limit=4.00,line
  width=4.560pt] (88.8500,8.5122) -- (58.8500,38.5122);

\end{tikzpicture}} & & &&
  ${}_*\mu_\XX^{\wedge}$& \scalebox{0.2}{\definecolor{c090000}{RGB}{9,0,0}
\definecolor{cff8080}{RGB}{255,128,128}

\begin{tikzpicture}[y=0.80pt, x=0.8pt,yscale=-1, inner sep=0pt, outer sep=0pt]
\path[cm={{1.0,0.0,0.0,-1.0,(248.85,-104.83783)}},draw=c090000,fill=cff8080,line
  cap=rect,miter limit=4.00,line width=4.560pt] (-130.0000,-93.3500) .. controls
  (-130.0000,-76.7815) and (-143.4315,-63.3500) .. (-160.0000,-63.3500) ..
  controls (-176.5685,-63.3500) and (-190.0000,-76.7815) .. (-190.0000,-93.3500)
  .. controls (-190.0000,-109.9185) and (-176.5685,-123.3500) ..
  (-160.0000,-123.3500) .. controls (-143.4315,-123.3500) and
  (-130.0000,-109.9185) .. (-130.0000,-93.3500) -- cycle;

\path[draw=black,line join=miter,line cap=butt,miter limit=4.00,line
  width=4.560pt] (108.8500,8.5122) -- (138.8500,38.5122);

\end{tikzpicture}}  \\
${}_*\mu_\XX^{\vee}$ & 
  \scalebox{0.2}{\definecolor{c090000}{RGB}{9,0,0}
\definecolor{cff8080}{RGB}{255,128,128}

\begin{tikzpicture}[y=0.80pt, x=0.8pt,yscale=-1, inner sep=0pt, outer sep=0pt]
\path[cm={{-1.0,0.0,0.0,1.0,(-51.15,101.86218)}},draw=c090000,fill=cff8080,line
  cap=rect,miter limit=4.00,line width=4.560pt]
  (-130.0000,-93.3500)arc(0.000:180.000:30.000)arc(-180.000:0.000:30.000) --
  cycle;

\path[draw=black,line join=miter,line cap=butt,miter limit=4.00,line
  width=4.560pt] (88.8500,-11.4878) -- (58.8500,-41.4878);

\end{tikzpicture}} & & &&
    ${}_*\mu_\XX^{^{/\!/}}$&   \scalebox{0.2}{\definecolor{c090000}{RGB}{9,0,0}
\definecolor{cff8080}{RGB}{255,128,128}

\begin{tikzpicture}[y=0.80pt, x=0.8pt,yscale=-1, inner sep=0pt, outer sep=0pt]
\path[shift={(248.85,101.86218)},draw=c090000,fill=cff8080,line cap=rect,miter
  limit=4.00,line width=4.560pt] (-130.0000,-93.3500) .. controls
  (-130.0000,-76.7815) and (-143.4315,-63.3500) .. (-160.0000,-63.3500) ..
  controls (-176.5685,-63.3500) and (-190.0000,-76.7815) .. (-190.0000,-93.3500)
  .. controls (-190.0000,-109.9185) and (-176.5685,-123.3500) ..
  (-160.0000,-123.3500) .. controls (-143.4315,-123.3500) and
  (-130.0000,-109.9185) .. (-130.0000,-93.3500) -- cycle;

\path[draw=black,line join=miter,line cap=butt,miter limit=4.00,line
  width=4.560pt] (108.8500,-11.4878) -- (138.8500,-41.4878);

\end{tikzpicture}} \\  
\end{tabular} 
 \end{center} 
 This is not always the case as we see in Example \ref{ex}. 
 \end{Remark}

The \textbf{parametrized homology} of $\XX$ is the collection of $\textrm{Dgm}_j^{\backslash \! \backslash}(\XX)$, $\textrm{Dgm}_j^{\vee}(\XX)$, $\textrm{Dgm}_j^{\wedge}(\XX)$, and $\textrm{Dgm}_j^{/ \! /}(\XX)$ over all $j$. 
Sometimes it is more convenient to use reduced instead of standard homology groups. We denote the four measures with respect to $\widetilde{\operatorname{\rm H}}_j$ by ${}_j\widetilde{\mu}_\XX^{{ \backslash\! \backslash}}$, 
${}_j\widetilde{\mu}_\XX^{\vee}$, ${}_j\widetilde{\mu}_\XX^{\wedge}$, and ${}_j\widetilde{\mu}_\XX^{^{/\!/}}$. The corresponding diagrams are $\operatorname{\rm \widetilde{D}gm}_j^{\backslash \! \backslash}(\XX)$, $\operatorname{\rm \widetilde{D}gm}_j^{\vee}(\XX)$, $\operatorname{\rm \widetilde{D}gm}_j^{\wedge}(\XX)$, and $\operatorname{\rm \widetilde{D}gm}_j^{/ \! /}(\XX)$. The \textbf{reduced parametrized homology} is the collection of these diagrams over all $j$.

We say that $\XX$ has a well-defined parametrized homology or a reduced parametrized homology when the four quantities defined above are finite r-measures. 
This is not always the case. Situations where $\XX = (X, p)$ has a well-defined parametrized homology include pairs when \cite{levelsetmeasure}:
\begin{itemize}
\item
 $X$ is a compact manifold with a boundary and $p$ is Morse, or more generally when $(X, p)$ is of Morse-type;
  \item
$X$ is a finite simplicial complex and $p$ is a piecewise-linear map;
\item
$X$ is a semialgebraic subset of $\RR^n$ and $p$ is the projection onto the $n$-th coordinate.
 \end{itemize}
\begin{Remark} 
In the case when $(X, p)$ is of Morse-type, leaving out the decorations on the points in the parametrized homology yields the levelset zigzag persistence diagram \cite{Zigzagpersistenthomologyandreal}.
\end{Remark}
 
For the four quantities to be r-measures they must be finite and additive with respect to horizontal and vertical splitting. The next two paragraphs summarize the proofs in \cite{levelsetmeasure}.

To show finiteness, let $R = [a, b]\times [c, d]$ be a rectangle in the half plane. Then for all $j$
$$
{}_j\mu^*(R) \leq \langle \scalebox{0.065}{\begin{tikzpicture}[y=0.80pt, x=0.8pt,yscale=-1, inner sep=0pt, outer sep=0pt]
\begin{scope}[shift={(-378.0,79.625)}]
  \path[cm={{-1.0,0.0,0.0,1.0,(480.0,0.0)}},draw=black,fill=black,line
    cap=rect,miter limit=4.00,fill opacity=0.000,line width=3.200pt]
    (100.0000,122.3622) .. controls (100.0000,138.9307) and (86.5685,152.3622) ..
    (70.0000,152.3622) .. controls (53.4315,152.3622) and (40.0000,138.9307) ..
    (40.0000,122.3622) .. controls (40.0000,105.7936) and (53.4315,92.3622) ..
    (70.0000,92.3622) .. controls (86.5685,92.3622) and (100.0000,105.7936) ..
    (100.0000,122.3622) -- cycle;

  \path[cm={{-1.0,0.0,0.0,1.0,(480.0,0.0)}},draw=black,fill=black,line
    cap=rect,miter limit=4.00,fill opacity=0.000,line width=3.200pt]
    (100.0000,122.3622) .. controls (100.0000,138.9307) and (86.5685,152.3622) ..
    (70.0000,152.3622) .. controls (53.4315,152.3622) and (40.0000,138.9307) ..
    (40.0000,122.3622) .. controls (40.0000,105.7936) and (53.4315,92.3622) ..
    (70.0000,92.3622) .. controls (86.5685,92.3622) and (100.0000,105.7936) ..
    (100.0000,122.3622) -- cycle;

  \path[cm={{-1.0,0.0,0.0,1.0,(480.0,0.0)}},draw=black,fill=black,line
    cap=rect,miter limit=4.00,fill opacity=0.000,line width=3.200pt]
    (100.0000,122.3622) .. controls (100.0000,138.9307) and (86.5685,152.3622) ..
    (70.0000,152.3622) .. controls (53.4315,152.3622) and (40.0000,138.9307) ..
    (40.0000,122.3622) .. controls (40.0000,105.7936) and (53.4315,92.3622) ..
    (70.0000,92.3622) .. controls (86.5685,92.3622) and (100.0000,105.7936) ..
    (100.0000,122.3622) -- cycle;

  \path[cm={{-1.0,0.0,0.0,1.0,(480.0,0.0)}},draw=black,fill=black,line
    cap=rect,miter limit=4.00,fill opacity=0.000,line width=3.200pt]
    (100.0000,122.3622) .. controls (100.0000,138.9307) and (86.5685,152.3622) ..
    (70.0000,152.3622) .. controls (53.4315,152.3622) and (40.0000,138.9307) ..
    (40.0000,122.3622) .. controls (40.0000,105.7936) and (53.4315,92.3622) ..
    (70.0000,92.3622) .. controls (86.5685,92.3622) and (100.0000,105.7936) ..
    (100.0000,122.3622) -- cycle;

  \path[shift={(340.0,0)},draw=black,fill=black,line cap=rect,miter
    limit=4.00,fill opacity=0.000,line width=3.200pt] (100.0000,122.3622) ..
    controls (100.0000,138.9307) and (86.5685,152.3622) .. (70.0000,152.3622) ..
    controls (53.4315,152.3622) and (40.0000,138.9307) .. (40.0000,122.3622) ..
    controls (40.0000,105.7936) and (53.4315,92.3622) .. (70.0000,92.3622) ..
    controls (86.5685,92.3622) and (100.0000,105.7936) .. (100.0000,122.3622) --
    cycle;

  \path[draw=black,line join=miter,line cap=butt,miter limit=4.00,draw
    opacity=0.551,line width=3.200pt] (430.0000,102.3622) -- (560.0000,-27.6378)
    -- (560.0000,-27.6378);

  \path[shift={(340.0,0)},draw=black,fill=black,line cap=rect,miter
    limit=4.00,fill opacity=0.000,line width=3.200pt] (100.0000,122.3622) ..
    controls (100.0000,138.9307) and (86.5685,152.3622) .. (70.0000,152.3622) ..
    controls (53.4315,152.3622) and (40.0000,138.9307) .. (40.0000,122.3622) ..
    controls (40.0000,105.7936) and (53.4315,92.3622) .. (70.0000,92.3622) ..
    controls (86.5685,92.3622) and (100.0000,105.7936) .. (100.0000,122.3622) --
    cycle;

  \path[draw=black,line join=miter,line cap=butt,miter limit=4.00,draw
    opacity=0.551,line width=3.200pt] (430.0000,102.3622) -- (560.0000,-27.6378)
    -- (560.0000,-27.6378);

  \path[shift={(340.0,0)},draw=black,fill=black,line cap=rect,miter
    limit=4.00,fill opacity=0.000,line width=3.200pt] (100.0000,122.3622) ..
    controls (100.0000,138.9307) and (86.5685,152.3622) .. (70.0000,152.3622) ..
    controls (53.4315,152.3622) and (40.0000,138.9307) .. (40.0000,122.3622) ..
    controls (40.0000,105.7936) and (53.4315,92.3622) .. (70.0000,92.3622) ..
    controls (86.5685,92.3622) and (100.0000,105.7936) .. (100.0000,122.3622) --
    cycle;

  \path[draw=black,line join=miter,line cap=butt,miter limit=4.00,draw
    opacity=0.551,line width=3.200pt] (430.0000,102.3622) -- (560.0000,-27.6378)
    -- (560.0000,-27.6378);

  \path[shift={(340.0,0)},draw=black,fill=black,line cap=rect,miter
    limit=4.00,line width=3.200pt] (100.0000,122.3622) .. controls
    (100.0000,138.9307) and (86.5685,152.3622) .. (70.0000,152.3622) .. controls
    (53.4315,152.3622) and (40.0000,138.9307) .. (40.0000,122.3622) .. controls
    (40.0000,105.7936) and (53.4315,92.3622) .. (70.0000,92.3622) .. controls
    (86.5685,92.3622) and (100.0000,105.7936) .. (100.0000,122.3622) -- cycle;

  \path[draw=black,line join=miter,line cap=butt,miter limit=4.00,line
    width=22.400pt] (430.0000,102.3622) -- (560.0000,-27.6378) --
    (560.0000,-27.6378);

  \path[cm={{-1.0,0.0,0.0,1.0,(650.0,-170.0)}},draw=black,fill=black,line
    cap=rect,miter limit=4.00,fill opacity=0.000,line width=3.200pt]
    (100.0000,122.3622) .. controls (100.0000,138.9307) and (86.5685,152.3622) ..
    (70.0000,152.3622) .. controls (53.4315,152.3622) and (40.0000,138.9307) ..
    (40.0000,122.3622) .. controls (40.0000,105.7936) and (53.4315,92.3622) ..
    (70.0000,92.3622) .. controls (86.5685,92.3622) and (100.0000,105.7936) ..
    (100.0000,122.3622) -- cycle;

  \path[draw=black,line join=miter,line cap=butt,miter limit=4.00,line
    width=22.400pt] (730.0000,102.3622) -- (600.0000,-27.6378) --
    (600.0000,-27.6378);

  \path[cm={{-1.0,0.0,0.0,1.0,(650.0,-170.0)}},draw=black,fill=black,line
    cap=rect,miter limit=4.00,line width=3.200pt] (100.0000,122.3622) .. controls
    (100.0000,138.9307) and (86.5685,152.3622) .. (70.0000,152.3622) .. controls
    (53.4315,152.3622) and (40.0000,138.9307) .. (40.0000,122.3622) .. controls
    (40.0000,105.7936) and (53.4315,92.3622) .. (70.0000,92.3622) .. controls
    (86.5685,92.3622) and (100.0000,105.7936) .. (100.0000,122.3622) -- cycle;

  \path[shift={(680.0,0)},draw=black,fill=black,line cap=rect,miter
    limit=4.00,line width=3.200pt] (100.0000,122.3622) .. controls
    (100.0000,138.9307) and (86.5685,152.3622) .. (70.0000,152.3622) .. controls
    (53.4315,152.3622) and (40.0000,138.9307) .. (40.0000,122.3622) .. controls
    (40.0000,105.7936) and (53.4315,92.3622) .. (70.0000,92.3622) .. controls
    (86.5685,92.3622) and (100.0000,105.7936) .. (100.0000,122.3622) -- cycle;

\end{scope}

\end{tikzpicture}} , \operatorname{\rm H}_j(\XX_{\{a, b, c, d\}}) \rangle = \textrm{dim}(\textrm{Im}(\operatorname{\rm H}_j(X_b) \to \operatorname{\rm H}_j(X_b^c))\cap \textrm{Im}(\operatorname{\rm H}_j(X_c) \to \operatorname{\rm H}_j(X_b^c) )).
$$
Bendich et al.\!\cite{bendich} show that $\textrm{Im}(\operatorname{\rm H}_j(X_b) \to \operatorname{\rm H}_j(X_b^c))\cap \textrm{Im}(\operatorname{\rm H}_j(X_c) \to \operatorname{\rm H}_j(X_b^c) )$ is the well group of $(X^{-1}([b, c]), p|_{X^{-1}([b, c])})$. Since its dimension is finite in the above situations, ${}_j\mu^*$ are finite.

The proof of additivity requires the Mayer-Vietoris principle for $X_a^b = X_a^p \cup X_p^b$ whenever $a<p<b$. This is automatically satisfied for a compact manifold and Morse function, as well as for a finite simplicial complex and piecewise-linear map. However, generally we have to restrict  $\XX$ to those whose level sets are embedded in a certain way. For example, this holds in cases where $X_b$ is a neighborhood deformation retract of $X_p^b$ or  of $X_b^q$  for all $b
\in \RR$ and some $p<b<q$. This restriction can be circumvented using a homology theory that satisfies the strong excision property. As a result, the four quantities are measures even in the case where $X$ is a compact simplicial complex and $p$ is a continuous function. 

\begin{Example}\label{example}
Consider the surface $X$ in Figure~\ref{primer}.  Since the projection $p$ onto the horizontal axis is Morse, $\XX$ has a well-defined parametrized homology. 
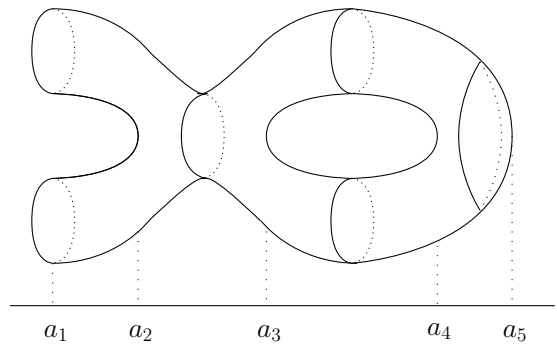
\begin{figure}[h!]
\begin{center}
\scalebox{0.5}{ \begin{tikzpicture}[y=0.80pt, x=0.8pt,yscale=-1, inner sep=0pt, outer sep=0pt]
\path[draw=black,line join=miter,line cap=butt,line width=0.800pt]
  (39.5122,1012.0370) -- (549.5122,1012.0370) -- (549.5122,1012.0370);

\path[draw=black,line join=miter,line cap=butt,line width=0.800pt]
  (79.5122,732.0370) .. controls (79.5122,732.0370) and (60.0000,732.0370) ..
  (60.0000,772.0370) .. controls (60.0000,812.0370) and (79.5122,812.0370) ..
  (79.5122,812.0370) -- (79.5122,812.0370);

\path[draw=black,dash pattern=on 0.80pt off 3.20pt,line join=miter,line
  cap=butt,miter limit=4.00,line width=0.800pt] (79.5122,732.0370) .. controls
  (79.5122,732.0370) and (100.0000,732.0370) .. (100.0000,772.0370) .. controls
  (100.0000,812.0370) and (79.5122,812.0370) .. (79.5122,812.0370) --
  (79.5122,812.0370);

\path[draw=black,line join=miter,line cap=butt,line width=0.800pt]
  (79.5122,892.0370) .. controls (79.5122,892.0370) and (60.0000,892.0370) ..
  (60.0000,932.0370) .. controls (60.0000,972.0370) and (79.5122,972.0370) ..
  (79.5122,972.0370) -- (79.5122,972.0370);

\path[draw=black,dash pattern=on 0.80pt off 3.20pt,line join=miter,line
  cap=butt,miter limit=4.00,line width=0.800pt] (79.5122,892.0370) .. controls
  (79.5122,892.0370) and (100.0000,892.0370) .. (100.0000,932.0370) .. controls
  (100.0000,972.0370) and (79.5122,972.0370) .. (79.5122,972.0370) --
  (79.5122,972.0370);

\path[draw=black,line join=miter,line cap=butt,line width=0.800pt]
  (79.5122,812.0370) .. controls (79.5122,812.0370) and (159.5122,812.0370) ..
  (159.5122,852.0370) .. controls (159.5122,892.0370) and (79.5122,892.0370) ..
  (79.5122,892.0370) -- (79.5122,892.0370) -- (79.5122,892.0370);

\path[draw=black,line join=miter,line cap=butt,line width=0.800pt]
  (224.5122,812.0370) .. controls (224.5122,812.0370) and (200.0000,812.0370) ..
  (200.0000,852.0370) .. controls (200.0000,892.0370) and (224.5122,892.0370) ..
  (224.5122,892.0370) -- (224.5122,892.0370);

\path[draw=black,line join=miter,line cap=butt,line width=0.800pt]
  (219.5122,812.0370);

\path[draw=black,line join=miter,line cap=butt,line width=0.800pt]
  (79.5122,812.0370) .. controls (79.5122,812.0370) and (159.5122,812.0370) ..
  (159.5122,852.0370) .. controls (159.5122,892.0370) and (79.5122,892.0370) ..
  (79.5122,892.0370) -- (79.5122,892.0370) -- (79.5122,892.0370);

\path[draw=black,line join=miter,line cap=butt,line width=0.800pt]
  (359.5122,812.0370) .. controls (359.5122,812.0370) and (279.5122,812.0370) ..
  (279.5122,852.0370) .. controls (279.5122,892.0370) and (359.5122,892.0370) ..
  (359.5122,892.0370) -- (359.5122,892.0370) -- (359.5122,892.0370);

\path[draw=black,line join=miter,line cap=butt,line width=0.800pt]
  (359.5122,812.0370) .. controls (359.5122,812.0370) and (439.5122,812.0370) ..
  (439.5122,852.0370) .. controls (439.5122,892.0370) and (359.5122,892.0370) ..
  (359.5122,892.0370) -- (359.5122,892.0370) -- (359.5122,892.0370);

\path[draw=black,line join=miter,line cap=butt,line width=0.800pt]
  (359.5122,972.0370) .. controls (359.5122,972.0370) and (509.5122,962.0370) ..
  (509.5122,852.0370) .. controls (509.5122,742.0370) and (359.5122,732.0370) ..
  (359.5122,732.0370) -- (359.5122,732.0370);

\path[draw=black,line join=miter,line cap=butt,line width=0.800pt]
  (359.5122,732.0370) .. controls (359.5122,732.0370) and (340.0000,732.0370) ..
  (340.0000,772.0370) .. controls (340.0000,812.0370) and (359.5122,812.0370) ..
  (359.5122,812.0370) -- (359.5122,812.0370);

\path[draw=black,line join=miter,line cap=butt,line width=0.800pt]
  (359.5122,892.0370) .. controls (359.5122,892.0370) and (340.0000,892.0370) ..
  (340.0000,932.0370) .. controls (340.0000,972.0370) and (359.5122,972.0370) ..
  (359.5122,972.0370) -- (359.5122,972.0370);

\path[draw=black,dash pattern=on 0.80pt off 3.20pt,line join=miter,line
  cap=butt,miter limit=4.00,line width=0.800pt] (219.5122,812.0370) .. controls
  (219.5122,812.0370) and (240.0000,812.0370) .. (240.0000,852.0370) .. controls
  (240.0000,892.0370) and (219.5122,892.0370) .. (219.5122,892.0370) --
  (219.5122,892.0370);

\path[draw=black,dash pattern=on 0.80pt off 3.20pt,line join=miter,line
  cap=butt,miter limit=4.00,line width=0.800pt] (359.5122,732.0370) .. controls
  (359.5122,732.0370) and (380.0000,732.0370) .. (380.0000,772.0370) .. controls
  (380.0000,812.0370) and (359.5122,812.0370) .. (359.5122,812.0370) --
  (359.5122,812.0370);

\path[draw=black,dash pattern=on 0.80pt off 3.20pt,line join=miter,line
  cap=butt,miter limit=4.00,line width=0.800pt] (359.5122,892.0370) .. controls
  (359.5122,892.0370) and (380.0000,892.0370) .. (380.0000,932.0370) .. controls
  (380.0000,972.0370) and (359.5122,972.0370) .. (359.5122,972.0370) --
  (359.5122,972.0370);

\path[draw=black,line join=miter,line cap=butt,line width=0.800pt]
  (79.5122,732.0370);

\path[draw=black,line join=miter,line cap=butt,line width=0.800pt]
  (169.5122,772.0370) .. controls (169.5122,772.0370) and (209.5122,812.0370) ..
  (219.5122,812.0370) .. controls (229.5122,812.0370) and (219.5122,812.0370) ..
  (219.5122,812.0370);

\path[draw=black,line join=miter,line cap=butt,line width=0.800pt]
  (149.5122,752.0370);

\path[draw=black,line join=miter,line cap=butt,line width=0.800pt]
  (169.5122,772.0370) .. controls (169.5122,772.0370) and (139.5122,732.0370) ..
  (79.5122,732.0370);

\path[draw=black,line join=miter,line cap=butt,line width=0.800pt]
  (169.5122,932.0370) .. controls (169.5122,932.0370) and (139.5122,972.0370) ..
  (79.5122,972.0370);

\path[draw=black,line join=miter,line cap=butt,line width=0.800pt]
  (269.5122,772.0370) .. controls (269.5122,772.0370) and (299.5122,732.0370) ..
  (359.5122,732.0370);

\path[draw=black,line join=miter,line cap=butt,line width=0.800pt]
  (274.5122,932.0370) .. controls (274.5122,932.0370) and (304.5122,972.0370) ..
  (364.5122,972.0370);

\path[draw=black,line join=miter,line cap=butt,line width=0.800pt]
  (269.5122,772.0370) .. controls (269.5122,772.0370) and (229.5122,812.0370) ..
  (219.5122,812.0370) .. controls (209.5122,812.0370) and (219.5122,812.0370) ..
  (219.5122,812.0370);

\path[draw=black,line join=miter,line cap=butt,line width=0.800pt]
  (169.5122,932.0370) .. controls (169.5122,932.0370) and (209.5122,892.0370) ..
  (219.5122,892.0370) .. controls (229.5122,892.0370) and (219.5122,892.0370) ..
  (219.5122,892.0370);

\path[draw=black,line join=miter,line cap=butt,line width=0.800pt]
  (273.9566,932.0370) .. controls (273.9566,932.0370) and (233.9566,892.0370) ..
  (223.9566,892.0370) .. controls (213.9566,892.0370) and (223.9566,892.0370) ..
  (223.9566,892.0370);

\path[draw=black,line join=miter,line cap=butt,line width=0.800pt]
  (479.5122,782.0370) .. controls (479.5122,782.0370) and (459.5122,812.0370) ..
  (459.5122,852.0370) .. controls (459.5122,892.0370) and (479.5122,922.0370) ..
  (479.5122,922.0370);

\path[draw=black,dash pattern=on 0.80pt off 4.80pt,line join=miter,line
  cap=butt,miter limit=4.00,line width=0.800pt] (479.5122,782.0370) .. controls
  (479.5122,782.0370) and (499.5122,812.0370) .. (499.5122,852.0370) .. controls
  (499.5122,892.0370) and (479.5122,922.0370) .. (479.5122,922.0370);

\path[draw=black,dash pattern=on 0.80pt off 6.40pt,line join=miter,line
  cap=butt,miter limit=4.00,line width=0.800pt] (79.5122,972.0370) --
  (79.5122,1012.0370);

\path[draw=black,dash pattern=on 0.80pt off 6.40pt,line join=miter,line
  cap=butt,miter limit=4.00,line width=0.800pt] (159.5122,942.0370) --
  (159.5122,1012.0370);

\path[draw=black,dash pattern=on 0.80pt off 6.40pt,line join=miter,line
  cap=butt,miter limit=4.00,line width=0.800pt] (279.5122,942.0370) --
  (279.5122,1012.0370) -- (279.5122,1012.0370);

\path[draw=black,dash pattern=on 0.80pt off 6.40pt,line join=miter,line
  cap=butt,miter limit=4.00,line width=0.800pt] (439.5000,955.3789) --
  (439.5000,1015.3789);

\path[draw=black,dash pattern=on 0.80pt off 6.40pt,line join=miter,line
  cap=butt,miter limit=4.00,line width=0.800pt] (509.5122,852.0370) --
  (509.5122,1012.0370) -- (509.5122,1012.0370);

\path[fill=black] (71.304062,1030.4875) node[below right] (text5404) {{\LARGE $a_1$}};

\path[fill=black] (150.67563,1030.4875) node[below right] (text3130) {{\LARGE $a_2$}};

\path[fill=black] (270.99936,1030.3961) node[below right] (text3138) {{\LARGE $a_3$}};

\path[fill=black] (430.12231,1028.5043) node[below right] (text3146) {{\LARGE $a_4$}};

\path[fill=black] (500.86606,1030.3961) node[below right] (text3154) {{\LARGE $a_5$}};

\end{tikzpicture} } 
\end{center}
\caption{\label{primer} Morse function on a 2-manifold with boundary.}
\end{figure}

To determine the diagrams belonging to each the four measures, we compute the multiplicities of the decorated points. When $p$ or $q$ is a regular point, the multiplicity of $(p^*, q^*)$ is $0$ for the four measures defined above. The only situations we have left to compute are when $p$ and $q$ are critical points. For example, we can now calculate the multiplicities of $(a_1^-, a_2^-)$ with respect to the four measures. Pick $\epsilon>0$ such that $a_1-\epsilon < a_1 <a_2-\epsilon < a_2$. We have
$$
\operatorname{\rm H}_0(\XX_{a_1-\epsilon, a_1, a_2-\epsilon, a_2}) \cong \scalebox{0.065}{\begin{tikzpicture}[y=0.80pt, x=0.8pt,yscale=-1, inner sep=0pt, outer sep=0pt]
\path[draw=black,fill=black,line cap=rect,miter limit=4.00,fill
  opacity=0.000,line width=3.200pt]
  (100.0000,122.3622)arc(0.000:180.000:30.000)arc(-180.000:0.000:30.000) --
  cycle;

\path[shift={(170.0,-170.0)},draw=black,fill=black,line cap=rect,miter
  limit=4.00,fill opacity=0.000,line width=3.200pt]
  (100.0000,122.3622)arc(0.000:180.000:30.000)arc(-180.000:0.000:30.000) --
  cycle;

\path[draw=black,fill=black,line cap=rect,miter limit=4.00,fill
  opacity=0.000,line width=3.200pt]
  (100.0000,122.3622)arc(0.000:180.000:30.000)arc(-180.000:0.000:30.000) --
  cycle;

\path[shift={(170.0,-170.0)},draw=black,fill=black,line cap=rect,miter
  limit=4.00,fill opacity=0.000,line width=3.200pt]
  (100.0000,122.3622)arc(0.000:180.000:30.000)arc(-180.000:0.000:30.000) --
  cycle;

\path[draw=black,fill=black,line cap=rect,miter limit=4.00,fill
  opacity=0.000,line width=3.200pt]
  (100.0000,122.3622)arc(0.000:180.000:30.000)arc(-180.000:0.000:30.000) --
  cycle;

\path[draw=black,line join=miter,line cap=butt,miter limit=4.00,draw
  opacity=0.316,line width=20.000pt] (90.0000,102.3622) -- (220.0000,-27.6378)
  -- (220.0000,-27.6378);

\path[shift={(170.0,-170.0)},draw=black,fill=black,line cap=rect,miter
  limit=4.00,fill opacity=0.000,line width=3.200pt]
  (100.0000,122.3622)arc(0.000:180.000:30.000)arc(-180.000:0.000:30.000) --
  cycle;

\path[draw=black,fill=black,line cap=rect,miter limit=4.00,fill
  opacity=0.000,line width=0.800pt]
  (100.0000,122.3622)arc(0.000:180.000:30.000)arc(-180.000:0.000:30.000) --
  cycle;

\path[shift={(170.0,-170.0)},draw=black,fill=black,line cap=rect,miter
  limit=4.00,fill opacity=0.000,line width=3.200pt]
  (100.0000,122.3622)arc(0.000:180.000:30.000)arc(-180.000:0.000:30.000) --
  cycle;

\path[cm={{-1.0,0.0,0.0,1.0,(480.0,0.0)}},draw=black,fill=black,line
  cap=rect,miter limit=4.00,fill opacity=0.000,line width=3.200pt]
  (100.0000,122.3622)arc(0.000:180.000:30.000)arc(-180.000:0.000:30.000) --
  cycle;

\path[draw=black,line join=miter,line cap=butt,miter limit=4.00,draw
  opacity=0.551,line width=3.200pt] (390.0000,102.3622) -- (260.0000,-27.6378)
  -- (260.0000,-27.6378);

\path[cm={{-1.0,0.0,0.0,1.0,(310.0,-170.0)}},draw=black,fill=black,line
  cap=rect,miter limit=4.00,fill opacity=0.000,line width=3.200pt]
  (100.0000,122.3622)arc(0.000:180.000:30.000)arc(-180.000:0.000:30.000) --
  cycle;

\path[cm={{-1.0,0.0,0.0,1.0,(480.0,0.0)}},draw=black,fill=black,line
  cap=rect,miter limit=4.00,fill opacity=0.000,line width=3.200pt]
  (100.0000,122.3622)arc(0.000:180.000:30.000)arc(-180.000:0.000:30.000) --
  cycle;

\path[draw=black,line join=miter,line cap=butt,miter limit=4.00,draw
  opacity=0.551,line width=3.200pt] (390.0000,102.3622) -- (260.0000,-27.6378)
  -- (260.0000,-27.6378);

\path[cm={{-1.0,0.0,0.0,1.0,(310.0,-170.0)}},draw=black,fill=black,line
  cap=rect,miter limit=4.00,fill opacity=0.000,line width=3.200pt]
  (100.0000,122.3622)arc(0.000:180.000:30.000)arc(-180.000:0.000:30.000) --
  cycle;

\path[cm={{-1.0,0.0,0.0,1.0,(480.0,0.0)}},draw=black,fill=black,line
  cap=rect,miter limit=4.00,fill opacity=0.000,line width=3.200pt]
  (100.0000,122.3622)arc(0.000:180.000:30.000)arc(-180.000:0.000:30.000) --
  cycle;

\path[draw=black,line join=miter,line cap=butt,miter limit=4.00,draw
  opacity=0.551,line width=3.200pt] (390.0000,102.3622) -- (260.0000,-27.6378)
  -- (260.0000,-27.6378);

\path[cm={{-1.0,0.0,0.0,1.0,(310.0,-170.0)}},draw=black,fill=black,line
  cap=rect,miter limit=4.00,fill opacity=0.000,line width=3.200pt]
  (100.0000,122.3622)arc(0.000:180.000:30.000)arc(-180.000:0.000:30.000) --
  cycle;

\path[cm={{-1.0,0.0,0.0,1.0,(480.0,0.0)}},draw=black,fill=black,line
  cap=rect,miter limit=4.00,fill opacity=0.000,line width=3.200pt]
  (100.0000,122.3622)arc(0.000:180.000:30.000)arc(-180.000:0.000:30.000) --
  cycle;

\path[draw=black,line join=miter,line cap=butt,miter limit=4.00,line
  width=22.400pt] (390.0000,102.3622) -- (260.0000,-27.6378) --
  (260.0000,-27.6378);

\path[cm={{-1.0,0.0,0.0,1.0,(310.0,-170.0)}},draw=black,fill=black,line
  cap=rect,miter limit=4.00,line width=3.200pt]
  (100.0000,122.3622)arc(0.000:180.000:30.000)arc(-180.000:0.000:30.000) --
  cycle;

\path[shift={(340.0,0)},draw=black,fill=black,line cap=rect,miter
  limit=4.00,fill opacity=0.000,line width=3.200pt]
  (100.0000,122.3622)arc(0.000:180.000:30.000)arc(-180.000:0.000:30.000) --
  cycle;

\path[draw=black,line join=miter,line cap=butt,miter limit=4.00,draw
  opacity=0.551,line width=3.200pt] (430.0000,102.3622) -- (560.0000,-27.6378)
  -- (560.0000,-27.6378);

\path[shift={(340.0,0)},draw=black,fill=black,line cap=rect,miter
  limit=4.00,fill opacity=0.000,line width=3.200pt]
  (100.0000,122.3622)arc(0.000:180.000:30.000)arc(-180.000:0.000:30.000) --
  cycle;

\path[draw=black,line join=miter,line cap=butt,miter limit=4.00,draw
  opacity=0.551,line width=3.200pt] (430.0000,102.3622) -- (560.0000,-27.6378)
  -- (560.0000,-27.6378);

\path[shift={(340.0,0)},draw=black,fill=black,line cap=rect,miter
  limit=4.00,fill opacity=0.000,line width=3.200pt]
  (100.0000,122.3622)arc(0.000:180.000:30.000)arc(-180.000:0.000:30.000) --
  cycle;

\path[draw=black,line join=miter,line cap=butt,miter limit=4.00,draw
  opacity=0.551,line width=3.200pt] (430.0000,102.3622) -- (560.0000,-27.6378)
  -- (560.0000,-27.6378);

\path[shift={(340.0,0)},draw=black,fill=black,line cap=rect,miter
  limit=4.00,line width=3.200pt]
  (100.0000,122.3622)arc(0.000:180.000:30.000)arc(-180.000:0.000:30.000) --
  cycle;

\path[draw=black,line join=miter,line cap=butt,miter limit=4.00,line
  width=22.400pt] (430.0000,102.3622) -- (560.0000,-27.6378) --
  (560.0000,-27.6378);

\path[cm={{-1.0,0.0,0.0,1.0,(650.0,-170.0)}},draw=black,fill=black,line
  cap=rect,miter limit=4.00,fill opacity=0.000,line width=3.200pt]
  (100.0000,122.3622)arc(0.000:180.000:30.000)arc(-180.000:0.000:30.000) --
  cycle;

\path[draw=black,line join=miter,line cap=butt,miter limit=4.00,line
  width=22.400pt] (730.0000,102.3622) -- (600.0000,-27.6378) --
  (600.0000,-27.6378);

\path[cm={{-1.0,0.0,0.0,1.0,(650.0,-170.0)}},draw=black,fill=black,line
  cap=rect,miter limit=4.00,line width=3.200pt]
  (100.0000,122.3622)arc(0.000:180.000:30.000)arc(-180.000:0.000:30.000) --
  cycle;

\path[shift={(680.0,0)},draw=black,fill=black,line cap=rect,miter
  limit=4.00,line width=3.200pt]
  (100.0000,122.3622)arc(0.000:180.000:30.000)arc(-180.000:0.000:30.000) --
  cycle;

\path[draw=black,line join=miter,line cap=butt,miter limit=4.00,draw
  opacity=0.314,line width=20.000pt] (770.0000,102.3622) -- (900.0000,-27.6378)
  -- (900.0000,-27.6378);

\path[cm={{-1.0,0.0,0.0,1.0,(990.0,-170.0)}},draw=black,fill=black,line
  cap=rect,miter limit=4.00,fill opacity=0.000,line width=3.200pt]
  (100.0000,122.3622)arc(0.000:180.000:30.000)arc(-180.000:0.000:30.000) --
  cycle;

\path[draw=black,line join=miter,line cap=butt,miter limit=4.00,draw
  opacity=0.314,line width=20.000pt] (1070.0000,102.3622) -- (940.0000,-27.6378)
  -- (940.0000,-27.6378);

\path[cm={{-1.0,0.0,0.0,1.0,(990.0,-170.0)}},draw=black,fill=black,line
  cap=rect,miter limit=4.00,fill opacity=0.000,line width=0.800pt]
  (100.0000,122.3622)arc(0.000:180.000:30.000)arc(-180.000:0.000:30.000) --
  cycle;

\path[shift={(1020.0,0)},draw=black,line cap=rect,miter limit=4.00,line
  width=0.800pt]
  (100.0000,122.3622)arc(0.000:180.000:30.000)arc(-180.000:0.000:30.000) --
  cycle;

\path[cm={{-1.0,0.0,0.0,1.0,(990.0,-170.0)}},draw=black,fill=black,line
  cap=rect,miter limit=4.00,fill opacity=0.000,line width=0.800pt]
  (100.0000,122.3622)arc(0.000:180.000:30.000)arc(-180.000:0.000:30.000) --
  cycle;

\end{tikzpicture}} \oplus  \scalebox{0.065}{\input{nn.tex}}  .
$$
The summand on the right is not registered by any of the measures, whereas the one on the left is detected by ${}_0\mu_\XX^{{ \backslash\! \backslash}}$. Since these values are the same for all $0<\epsilon< a_2 - a_1$, we have
$$
\begin{array}{lclcccccc}
m_{{}_0\mu_\XX^{{ \backslash\! \backslash}} }(a_1^-, a_2^-)&= & \displaystyle\lim_{\epsilon \to 0} {}_0\mu_\XX^{{ \backslash\! \backslash}}([a_1-\epsilon, a_1] \times [a_2 -\epsilon, a_2]) & = & 1,\\
m_{{}_0\mu_\XX^{\vee}}(a_1^-, a_2^-)& = & \displaystyle\lim_{\epsilon \to 0} {}_0\mu_\XX^{\vee}([a_1-\epsilon, a_1] \times [a_2-\epsilon, a_2])&=&0,\\
m_{{}_0\mu_\XX^{\wedge}} (a_1^-, a_2^-)&= & \displaystyle\lim_{\epsilon \to 0} {}_0\mu_\XX^{\wedge}([a_1-\epsilon, a_1] \times [a_2-\epsilon, a_2])&=&0,\\
m_{{}_0\mu_\XX^{^{/\!/}} }(a_1^-, a_2^-)&= &\displaystyle\lim_{\epsilon \to 0} {}_0\mu_\XX^{^{/\!/}}([a_1-\epsilon, a_1] \times [a_2-\epsilon, a_2])&=&0.
\end{array}
$$
This means that $(a_1^-, a_2^-)$ is a point in the decorated persistence diagram belonging to ${}_0\mu_\XX^{{ \backslash\! \backslash}}$ with a multiplicity of 1. The corresponding 0-homology cycle ceases to exist beyond $a_1$, and is killed at $a_2$.

We repeat this procedure to compute the other multiplicities.
The parametrized homology of $\XX$ is represented in Figure~\ref{Diagram}.
\begin{figure}[h!]
\begin{center}
\includegraphics[scale=0.3]{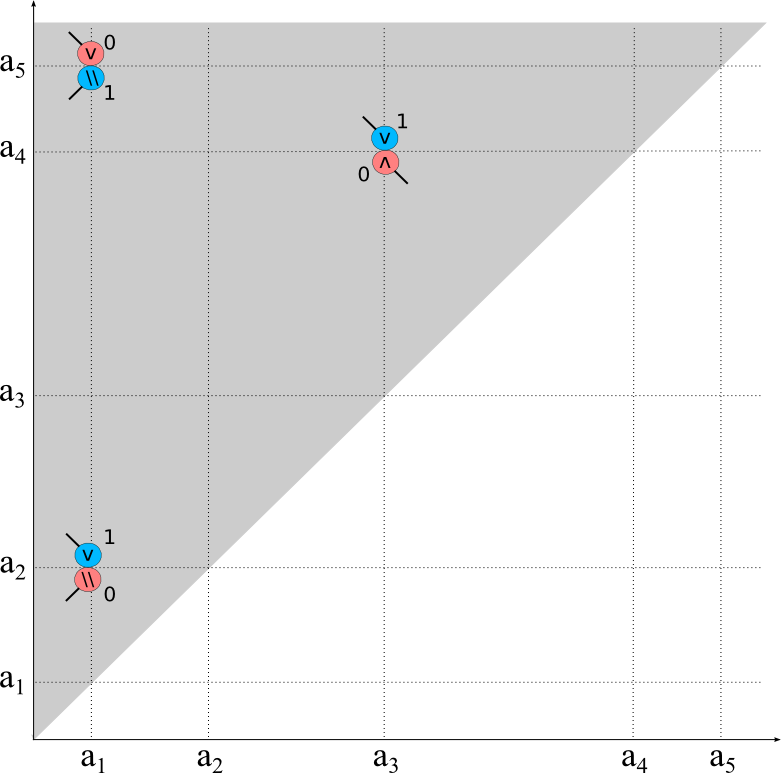}
 \end{center}
\caption{\label{Diagram}The parametrized homology of $\XX$.  The color of the point indicates the dimension, while the symbol designates to which of the  diagrams it belongs.}
\end{figure}

\end{Example}

\section{Alexander Duality for Parametrized Homology}
Alexander duality is a statement about the relationship between the cohomology groups of a locally contractible, compact subset of $\RR^n$ and the homology groups of the complement (for more on the classical statement, see Hatcher \cite{hatcher} and Spanier  \cite{spanier} ). 
\begin{Theorem}[Alexander Duality]
If $K$ is a locally contractible, compact subset of $\RR^n$, then for all $j = 0, \ldots,  n-1$,
$$
\widetilde{\operatorname{\rm H}}_{n-j-1}(\RR^n-K) \cong \operatorname{\rm H}^j(K).
$$
\end{Theorem}

The goal of this section is to extend this theorem to parametrized spaces. 
\begin{Theorem}\label{theorem}
Let  $X \subset \RR^n \times \RR$ with $n\geq 2$ be a compact set, let $Y = (\RR^n \times \RR) \setminus X$,
and let $p$ be the projection onto the second factor. We assume that the level sets $X_a$ for $a\in \RR$, and slices $X_a^b$ for $a<b$ are locally contractible. If $(X, p|_X)$ has a well-defined parametrized homology, then  the pair $(Y, p|_Y)$ has a well-defined reduced parametrized homology. Additionally, for $j  = 0, \ldots, n-1$:
$$
\begin{array}{lcr}
\operatorname{\rm \widetilde{D}gm}_{n-j-1}^{\backslash\! \backslash}(\YY) &= & \operatorname{\rm Dgm}_{j}^{/\!/}(\XX)  \\
\operatorname{\rm \widetilde{D}gm}_{n-j-1}^{\vee}(\YY) &= & \operatorname{\rm Dgm}_{j}^{\wedge}(\XX)  \\
\operatorname{\rm \widetilde{D}gm}_{n-j-1}^{\wedge}(\YY) &= & \operatorname{\rm Dgm}_{j}^{\vee}(\XX)  \\
\operatorname{\rm \widetilde{D}gm}_{n-j-1}^{/\!/}(\YY) &= & \operatorname{\rm Dgm}_{j}^{\backslash\! \backslash}(\XX)
\end{array}
$$
(recall that $X_a = p^{-1}(a) \cap X$ and $X_a^b =p^{-1}([a,b]) \cap X$).
\end{Theorem}
\begin{Remark}\label{rem}
From the proof we can deduce the following duality: if a $j$-dimensional homology cycle in $\XX$ is killed (ceases to exist) at endpoint $p$, then there is a corresponding $(n-j-1)$-dimensional homology cycle in $\YY$, which ceases to exist (is killed) beyond that same endpoint.
\end{Remark}
\begin{Remark}
The conditions of the theorem are satisfied for $(X, p|_X)$, where:
\begin{itemize}
\item
 $X$ is a compact submanifold of $\RR^n\times \RR$ (with or without boundary) and $p|_X$ is Morse; or, more generally, when $(X, p|_X)$ is of Morse-type and $X$ is compact with all the slices and level sets locally contractible; 
\item
 $X$ is a finite simplicial complex and $p|_X$ is a piecewise-linear map;
 \item
 $X$ is a semialgebraic subset of $\RR^n\times \RR$.
 \end{itemize}
The methods of this paper clear the way for a more general theorem using a homology functor with strong excision, such as the \v{C}ech cohomology.
\end{Remark}
Example \ref{34} gives an instance of our result not covered by the Land and Water theorem by Edelsbrunner and Kerber \cite{AlexanderDuality}:
\begin{Example}\label{34}
Let $S \subset \RR^3$ be the Alexander horned sphere (see \cite[Example 2B.2]{hatcher}). Let $X = S \times [-1,1] \subset \RR^3 \times \RR$, let $Y = \RR^3\times \RR \setminus X$, and let $p$ be the projection onto the second factor. Since $(X, p|_X)$ is of Morse-type, it has a well-defined  parametrized homology. The conditions of Theorem \ref{theorem} are satisfied, because $S$ is locally contractible and compact. 
\end{Example}

The proof of Theorem \ref{theorem} requires two lemmas.

\begin{Lemma}\label{exact}
Let  $X$, $Y$, and $p$ be as in the theorem. Consider the following diagram of vector spaces and maps:
\begin{center}
\begin{tikzpicture}[font=\small, xscale=1.2,yscale=1.2]
\draw (0,2) node(02){} (1,2) node(12){$\widetilde{\operatorname{\rm H}}_{n-j-1}(Y_{a}^{b})$} (2,2) node(22){} ;
\draw (-1,1) node(01){$\widetilde{\operatorname{\rm H}}_{n-j-1}(Y_a)$} (1,1) node(11){} (3,1) node(21){$\widetilde{\operatorname{\rm H}}_{n-j-1}(Y_b)$};
\draw (-1,0) node(00){$ \operatorname{\rm H}^{j}(X_a)$} (1,0) node(10){} (3,0) node(20){$\operatorname{\rm H}^{j}(X_b)$};
\draw (0,-1) node(0){} (1,-1) node(1){$\operatorname{\rm H}^{j}(X_a^b)$} (2,-1) node(2){} ;

\draw[->] (01) to node[above]{$i_a$} (12); 
\draw[->] (21) to node[above]{$i_b$} (12); 
\draw[->] (1) to node[below]{$i^a$} (00);
\draw[->] (1) to node[below]{$i^b$} (20); 
\draw[->] (00) to node[right]{$D_a$}(01);
\draw[->] (20) to node[right]{$D_b$} (21);
\end{tikzpicture}
\end{center}
Maps $i^a, i^b, i_a$, and $i_b$ are induced by the inclusions $X_a \hookrightarrow X_a^b$, $X_b \hookrightarrow X_a^b$,  $Y_a \hookrightarrow Y_a^b$, and $Y_b \hookrightarrow Y_a^b$. Isomorphisms $D_a$ and $D_b$, respectively, are Alexander duality isomorphisms in $\RR^n \times \{a\}$ and $\RR^n \times \{b\}$. Then $\textrm{Im}\, (D_ai^a \oplus D_b i^b) = \textrm{Ker}\, (i_a - i_b)$. 
\end{Lemma}
\begin{Remark}
This lemma holds even if we do not assume field coefficients. 
\end{Remark}

\begin{proof}  
We look at the long exact sequence for homology groups of the pair $(Y_a^b, Y_a \cup Y_b)$:
\begin{equation}
\cdots \to \operatorname{\rm H}_{n-j}(Y_a^b, Y_a \cup Y_b)\to \operatorname{\rm H}_{n-j-1}(Y_a \cup Y_b) \to\operatorname{\rm H}_{n-j-1}(Y_a^b)\to \operatorname{\rm H}_{n-j-1}(Y_a^b, Y_a \cup Y_b) \to \cdots
\end{equation}
We claim that there exists an isomorphism $\operatorname{\rm H}^{j}(X_a^b) \to \operatorname{\rm H}_{n-j}(Y_a^b, Y_a \cup Y_b)$ making the following diagram commute up to a sign
\begin{center}
\begin{tikzpicture}[xscale=5,yscale=1, font=\small]
\draw (1.5,-1) node(1){} (2,-1) node(2){$\operatorname{\rm H}^{j}(X_a^b)$} (3,-1) node(3){$\operatorname{\rm H}^{j}(X_a) \oplus \operatorname{\rm H}^{j}(X_b)$} (4,-1) node(4){} (4.5,-1) node(5){};
\draw (1.5,-2) node(91){} (2,-2) node(92){$\operatorname{\rm H}_{n-j}(Y_a^b, Y_a \cup Y_b)$} (3,-2) node(93){$\operatorname{\rm H}_{n-j-1}(Y_a)\oplus \operatorname{\rm H}_{n-j-1}(Y_b)$} (4,-2) node(94){$\operatorname{\rm H}_{n-j-1}(Y_a^b).$} (4.5,-2) node(95){};

\draw[->] (2) to node[above] {$i^a\oplus i^b$} (3); 
\draw[->](93) to node[above]{$i_a-i_b$} (94);
 
\draw[->] (92) -- (93);

\draw[->] (2) to node[left]{\rotatebox[origin=c]{90}{$\cong$}}(92);
\draw[->] (3) to node[right]{\rotatebox[origin=c]{90}{$\cong$}} (93);

\end{tikzpicture}
\end{center}
Granted that, the case $j \neq n-1$ follows immediately by virtue of exactness of the bottom line. We analyze the case when $j=n-1$ separately. 

Let $\operatorname{\rm H}_c^*$ denote cohomology with compact supports. According to \cite[Chapter 3, Problem 35]{hatcher}, the following diagram, where the horizontal
lines are long exact sequences of the corresponding pairs and the vertical
arrows are Poincar\'{e} duality isomorphisms, commutes up to a sign.
\begin{center}
\begin{tikzpicture}[xscale=5,yscale=1, font=\small]
\draw (1.5,-1) node(1){$\cdots$} (2,-1) node(2){$\operatorname{\rm H}^{j+1}_c(Y_a^b)$} (3,-1) node(3){$\operatorname{\rm H}^{j+1}_c(Y_a \cup Y_b)$} (4,-1) node(4){$ \operatorname{\rm H}^{j+2}_c(Y_a^b, Y_a \cup Y_b)$} (4.5,-1) node(5){$\cdots$};
\draw (1.5,-2) node(91){$\cdots$} (2,-2) node(92){$\operatorname{\rm H}_{n-j}(Y_a^b, Y_a \cup Y_b)$} (3,-2) node(93){$\operatorname{\rm H}_{n-j-1}(Y_a \cup Y_b)$} (4,-2) node(94){$\operatorname{\rm H}_{n-j-1}(Y_a^b)$} (4.5,-2) node(95){$\cdots$};

\draw[->] (1) --  (2); 
\draw[->] (2) -- (3); 
\draw[->] (3) -- (4);
\draw[->] (4) -- (5); 
 
\draw[->] (91) -- (92); 
\draw[->] (92) -- (93);
\draw[->] (93) -- (94); 
\draw[->] (94) -- (95);

\draw[->] (2)  to node[left]{\rotatebox[origin=c]{90}{$\cong$}}(92);
\draw[->] (3)  to node[left]{\rotatebox[origin=c]{90}{$\cong$}}(93);
\draw[->] (4)  to node[left]{\rotatebox[origin=c]{90}{$\cong$}} (94);

\end{tikzpicture}
\end{center}

Let $\{ U_i\}_i$ be a nested sequence of neighborhoods of $X_a^b$ in $\RR^n\times[a, b]$, such that $U_1$ retracts onto $X_a^b$ and $\cap_i U_i = X_a^b$. Such sequences exist since $X_a^b$ is compact and
locally contractible by Theorem A.7 of \cite{hatcher}. Further let $\{B_i\}$ be an increasing sequence
of closed balls centered at the origin and containing $X_a^b$ such that
$\cup_iB_i=\mathbb{R}^{n+1}$.  We may assume that
$\overline U_1\subset\mathrm{Int} B_1$ so that
$\overline U_i\cap\overline{B_i^C}=\emptyset$ for all $i$.
 Now $U_1 \cap \RR^n \times \{a\}$ is open in $\RR^n\times \{a\}$ in the subspace topology and contains $X_a$. Since $X_a$ is compact and locally contractible, we can find a neighborhood $U^a$ of $X_a$ which retracts onto $X_a$ (again using Theorem A.7 of \cite{hatcher}). Pick
a nested sequence of neighborhoods $U^a_i$ of $X_a$ such that
$U^a_i\subset U^a\cap U_i$ for each $i$ and $\cap_iU^a_i=X_a$. In a similar manner we obtain a system of neighborhoods for $X_b$. 

Let $A^C$ denote the complement of $A$ (where the ambient set is clear
from the context). By cofinality, we have
$$
{\small
\begin{array}{rclcl}
\operatorname{\rm H}^{j+1}_c(Y_a^b) &\cong& \colim \operatorname{\rm H}^{j+1}(Y_a^b, Y_a^b \setminus B_i \cap U_i^C), \\
\operatorname{\rm H}^{j+1}_c(Y_a \cup Y_b)& \cong& \colim \operatorname{\rm H}^{j+1}(Y_a, Y_a\setminus B_i \cap (U_i^a)^C)\oplus \operatorname{\rm H}^{j+1}(Y_b, Y_b \setminus B_i \cap (U_i^b)^C) \\
\end{array}}
$$
Moreover, the restriction $ \operatorname{\rm H}^{j+1}_c(Y_a^b)  \to \operatorname{\rm H}^{j+1}_c(Y_a \cup Y_b)$
is the the colimit of the corresponding morphisms.

Using the notation $(B_i^C)_a^b = B_i^C \cap \RR^n\times [a,b]$, $(B_i^C)_a = B_i^C \cap \RR^n\times\{a\}$, and $(B_i^C)_b = B_i^C \cap \RR^n\times\{b\}$, we rewrite the expressions
$$
{\small
\begin{array}{rcl}
\operatorname{\rm H}^{j+1}(Y_a^b, Y_a^b \setminus B_i \cap U_i^C) & = & \operatorname{\rm H}^{j+1}(Y_a^b, (B_i^C)_a^b \cup (U_i\setminus X_a^b)),\\
\operatorname{\rm H}^{j+1}(Y_a, Y_a\setminus B_i \cap (U_i^a)^C)&=& \operatorname{\rm H}^{j+1}(Y_a, (B_i^C)_a \cup (U_i^a\setminus X_a)) ,\\
\operatorname{\rm H}^{j+1}(Y_b, Y_b \setminus B_i \cap (U_i^b)^C)  & = &  \operatorname{\rm H}^{j+1}(Y_b, (B_i^C)_b \cup (U_i^b\setminus X_b)).\\
\end{array} }
$$
We have
$$
{\small
\begin{array}{rclcl}
\operatorname{\rm H}^{j+1}(Y_a^b, (B_i^C)_a^b \cup (U_i\setminus  X_a^b)) &\cong &  \operatorname{\rm H}^{j+1}(\RR^n\times [a,b], (B_i^C)_a^b \cup U_i) &  \cong & \widetilde{\operatorname{\rm H}}^{j}((B_i^C)_a^b \cup U_i),\\
\operatorname{\rm H}^{j+1}(Y_a, (B_i^C)_a \cup (U_i^a\setminus X_a))&\cong & \operatorname{\rm H}^{j+1}(\RR^n\times \{a\}, (B_i^C)_a \cup U_i^a) & \cong & \widetilde{\operatorname{\rm H}}^{j}((B_i^C)_a \cup U_i^a), \\
\operatorname{\rm H}^{j+1}(Y_b, (B_i^C)_b \cup (U_i^a\setminus X_b))&\cong & \operatorname{\rm H}^{j+1}(\RR^n\times \{b\}, (B_i^C)_b \cup U_i^b) &\cong & \widetilde{\operatorname{\rm H}}^{j}((B_i^C)_b \cup U_i^b) .\\
\end{array} }
$$
The left-hand isomorphisms follow from excision and the right-hand isomorphisms from the long exact sequence of a pair. 

By naturality of the above isomorphisms and by commutativity of the
Poincar\'{e} duality ladder, the following diagram commutes up to a sign.

 \begin{tikzpicture}[font=\small] 
\draw (0,1) node(2){$\colim \widetilde{\operatorname{\rm H}}^{j}((B_i^C)_a^b \cup U_i)$} (8,1) node(3){$\colim (\widetilde{\operatorname{\rm H}}^{j}((B_i^C)_a \cup U_i^a) \oplus \widetilde{\operatorname{\rm H}}^{j}((B_i^C)_b \cup U_i^b))$} (8,1) ;
\draw (14.1, 0.5) node(40){ {\normalsize (2)}};
\draw (0,0) node(92){$\operatorname{\rm H}_{n-j}(Y_a^b, Y_a \cup Y_b)$} (8,0) node(93){$\operatorname{\rm H}_{n-j-1}(Y_a \cup Y_b).$} (8,0)  ;

\draw[->] (2) -- (3); 
 
\draw[->] (92) -- (93);

\draw[->] (2)  to node[left]{\rotatebox[origin=c]{90}{$\cong$}} (92);
\draw[->] (3)  to node[right]{\rotatebox[origin=c]{90}{$\cong$}}(93);

\end{tikzpicture}

Arguing as Hatcher, in Theorem 3.44   \cite{hatcher}, we infer
\setcounter{equation}{2}
\begin{equation}
\colim \widetilde{\operatorname{\rm H}}^{j}((B_i^C)_a^b \cup U_i) \cong \colim ( \widetilde{\operatorname{\rm H}}^j(S^{n-1})) \oplus \colim( \operatorname{\rm H}^{j}(U_i))  \cong  \widetilde{\operatorname{\rm H}}^j(S^{n-1}) \oplus \operatorname{\rm H}^{j}(X_a^b),
\end{equation}
and similarly
\begin{equation}
\colim (\widetilde{\operatorname{\rm H}}^{j}((B_i^C)_a \cup U_i^a) \oplus \widetilde{\operatorname{\rm H}}^{j}((B_i^C)_b \cup U_i^b))\cong  \widetilde{\operatorname{\rm H}}^j(S^{n-1}) \oplus \operatorname{\rm H}^{j}(X_a) \oplus \widetilde{\operatorname{\rm H}}^j(S^{n-1}) \oplus \operatorname{\rm H}^{j}(X_b).
\end{equation}
If $j\neq n-1$, we have $\widetilde{\operatorname{\rm H}}^j(S^{n-1}) =0$. We insert (3) and (4)
into (2) and get the desired commutative square
\begin{center}
\begin{tikzpicture}[xscale=8,yscale=1, font=\small]
\draw (1.5,-1) node(1){} (2,-1) node(2){$\operatorname{\rm H}^{j}(X_a^b)$} (3,-1) node(3){$\operatorname{\rm H}^{j}(X_a) \oplus \operatorname{\rm H}^{j}(X_b)$} (4,-1) node(4){} (4.5,-1) node(5){};
\draw (1.5,-2) node(91){} (2,-2) node(92){$\operatorname{\rm H}_{n-j}(Y_a^b, Y_a \cup Y_b)$} (3,-2) node(93){$\operatorname{\rm H}_{n-j-1}(Y_a)\oplus \operatorname{\rm H}_{n-j-1}(Y_b).$} (4,-2) node(94){} (4.5,-2) node(95){};

\draw[->] (2) -- (3); 
 
\draw[->] (92) -- (93);

\draw[->] (2)  to node[left]{\rotatebox[origin=c]{90}{$\cong$}} (92);
\draw[->] (3)  to node[right]{\rotatebox[origin=c]{90}{$\cong$}} (93);

\end{tikzpicture}
\end{center}
This finishes the proof for $j\neq n-1$.

Now let $j=n-1$. In this case $\widetilde{\operatorname{\rm H}}^j(S^{n-1}) \cong \mathbf k $. Observe that 
$$
\operatorname{\rm H}_0(Y_a \cup Y_b) \cong \mathbf k \oplus \widetilde{\operatorname{\rm H}}_0 (Y_a) \oplus \mathbf k \oplus \widetilde{\operatorname{\rm H}}_0 (Y_b) \quad
\textrm{and}\quad
\operatorname{\rm H}_0(Y_a^b) \cong \mathbf k \oplus  \widetilde{\operatorname{\rm H}}_0 (Y_a^b).
$$
Taking this into account, inserting (3) and (4) into (2), and extending the bottom line by an extra term from (1), we get the following commutative diagram
\begin{center}
\begin{tikzpicture}[xscale=5,yscale=1, font=\small]
\draw (1.5,-1) node(1){} (2,-1) node(2){$\mathbf k \oplus \operatorname{\rm H}^{n-1}(X_a^b)$} (3,-1) node(3){$\mathbf k \oplus \operatorname{\rm H}^{n-1}(X_a) \oplus \mathbf k  \oplus \operatorname{\rm H}^{n-1}(X_b)$} (4,-1) node(4){} (4.5,-1) node(5){};
\draw (1.5,-2) node(91){} (2,-2) node(92){$\operatorname{\rm H}_{1}(Y_a^b, Y_a \cup Y_b)$} (3,-2) node(93){$\mathbf k \oplus  \widetilde{\operatorname{\rm H}}_{0}(Y_a)\oplus \mathbf k \oplus \widetilde{\operatorname{\rm H}}_{0}(Y_b)$} (4,-2) node(94){$\mathbf k \oplus \widetilde{\operatorname{\rm H}}_{0}(Y_a^b).$} (4.5,-2) node(95){};

\draw[->] (2) to node[above] {} (3); 
\draw[->](93) to node[above]{} (94);
 
\draw[->] (92) -- (93);

\draw[->] (2) to node[left]{\rotatebox[origin=c]{90}{$\cong$}}(92);
\draw[->] (3) to node[right]{\rotatebox[origin=c]{90}{$\cong$}} (93);

\end{tikzpicture}
\end{center}
Once again, the vertical maps are isomorphisms.  Since additional copies of $\mathbf k$ get mapped to corresponding additional
copies of $\mathbf k$, exactness of the diamond diagram from the statement of
the lemma now follows for $j=n-1$.
\end{proof}
We also need the Diamond Principle \cite{Zigzagpersistence} for the proof of theorem \ref{theorem}. Consider the following diagram of vector spaces and linear maps between them.
\begin{center}
\begin{tikzpicture}[font=\small, xscale=1.5,yscale=1.5]

\draw (0,0) node(00){$V_1$} ;
\draw (1,0) node(10){$\ldots$} ;
\draw (2,0) node(20){$V_{k-1}$} ;
\draw (3,1) node(31){$W_k$} ;
\draw (3,-1) node(301){$U_k$};
\draw (4,0) node(40){$V_{k+1}$};
\draw (5,0) node(50){$\ldots$};
\draw (6,0) node(60){$V_{n}$};

\draw[<->] (00) to node[above]{$p_1$} (10); 
\draw[<->] (10) to node[above]{$p_{k-2}$} (20);
\draw[->] (20) to node[left]{$f_{k-1}$} (31);
\draw[->] (40) to node[right]{$g_k$} (31); 
\draw[->] (301) to node[left]{$g_{k-1}$} (20);
\draw[->] (301) to node[right]{$f_k$} (40); 
\draw[<->] (40) to node[above]{$p_{k+1}$} (50); 
\draw[<->] (50) to node[above]{$p_{n-1}$} (60); 

\end{tikzpicture}
\end{center}

We say that the diamond in the center 
is exact if $\textrm{Im}(D_1) = \textrm{Ker}(D_2)$ in the following sequence
$$
\begin{tikzpicture}[font=\small, xscale=1,yscale=1]
\draw (0,0) node(00){$U_k$} ;
\draw (2,0) node(10){$V_{k-1}\oplus V_{k+1}$} ;
\draw (4,0) node(20){$W_k,$} ;
\draw[->] (00) to node[above]{$D_1$} (10);
\draw[->] (10) to node[above]{$D_2$} (20); 
\end{tikzpicture}
$$
where $D_1(u) = g_{k-1}(u) \oplus f_k(u)$, and $D_2(v \oplus v') = f_{k-1}(v)- g_k(v')$.

Let $\VV^+$ and $\VV^-$ denote the upper and lower zigzag modules.
$$\VV^+ =
\begin{tikzpicture}[font=\small, xscale=1,yscale=1]
\draw (0,0) node(00){$V_1$} ;
\draw (2,0) node(10){$\ldots$} ;
\draw (4,0) node(20){$V_{k-1}$} ;
\draw (6,0) node(31){$W_k$} ;
\draw (8,0) node(40){$V_{k+1}$};
\draw (10,0) node(50){$\ldots$};
\draw (12,0) node(60){$V_{n},$};

\draw[<->] (00) to node[above]{$p_1$} (10); 
\draw[<->] (10) to node[above]{$p_{k-2}$} (20);
\draw[->] (20) to node[above]{$f_{k-1}$} (31);
\draw[->] (40) to node[above]{$g_k$} (31); 
\draw[<->] (40) to node[above]{$p_{k+1}$} (50); 
\draw[<->] (50) to node[above]{$p_{n-1}$} (60); 

\end{tikzpicture}
$$
$$\VV^- =
\begin{tikzpicture}[font=\small, xscale=1,yscale=1]
\draw (0,0) node(00){$V_1$} ;
\draw (2,0) node(10){$\ldots$} ;
\draw (4,0) node(20){$V_{k-1}$} ;
\draw (6,0) node(31){$U_k$} ;
\draw (8,0) node(40){$V_{k+1}$};
\draw (10,0) node(50){$\ldots$};
\draw (12,0) node(60){$V_{n}.$};

\draw[<->] (00) to node[above]{$p_1$} (10); 
\draw[<->] (10) to node[above]{$p_{k-2}$} (20);
\draw[->] (31) to node[above]{$g_{k-1}$} (20);
\draw[->] (31) to node[above]{$f_k$} (40); 
\draw[<->] (40) to node[above]{$p_{k+1}$} (50); 
\draw[<->] (50) to node[above]{$p_{n-1}$} (60); 
\end{tikzpicture}
$$
We have the following relation between persistence diagrams of $\VV^+$ and $\VV^-$.
\begin{Lemma}[The Diamond Principle \cite{Zigzagpersistence}]
Given $\VV^+$ and $\VV^-$ as above, suppose that the middle diamond is exact. Then there is a partial bijection between the set of intervals that appear in the interval modules decomposition of $\VV^+$ and the set of intervals that appear in the interval modules decomposition of $\VV^-$. Intervals are matched according to the following rules:
\begin{itemize}
\item
Intervals of type $[k, k]$ are unmatched,
\item
 Type $[b,k]$ is matched with type $[b, k-1]$ and vice versa, for $b\leq k-1$,
 \item 
 Type $[k,d]$ is matched with type $[k+1,d]$ and vice versa, for $d\geq k+1$,
\item
  Type $[b, d]$ is matched with type $[b, d]$ in all other cases.
\end{itemize}
\begin{center}
\includegraphics[scale=0.12]{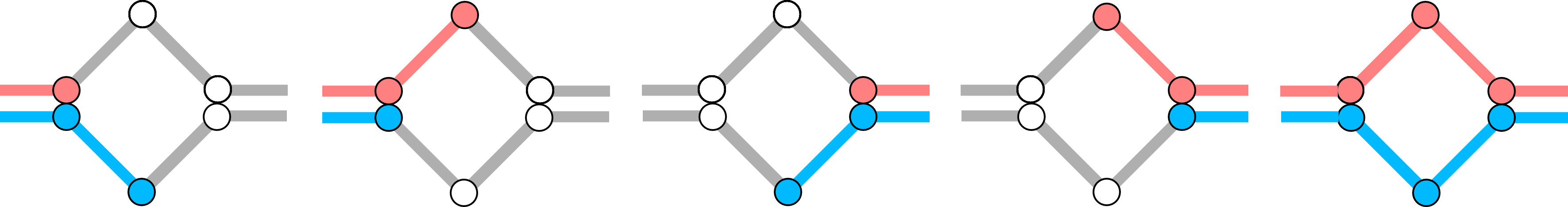}
\end{center}
\end{Lemma}

\begin{proof}[Proof of Theorem \ref{theorem}]
By assumption, $(X, p|_X)$ has a well-defined parametrized homology, so ${}_j\mu_\XX^{^{\backslash\! \backslash}}$, ${}_j\mu_\XX^{^{\vee}}$, ${}_j\mu_\XX^{^{\wedge}}$, and ${}_j\mu_\XX^{^{/\!/}}$ are finite r-measures for $j=0, \ldots, n-1$. Let $R = [a, b]\times [c, d]$ with $-\infty <a<b<c<d<\infty$. By Lemma~\ref{exact} the three diamonds in the following diagram are exact:
\begin{center}
\begin{tikzpicture}[font=\small, xscale=1.2,yscale=1.2]

\draw (-3,2) node(032){$\widetilde{\operatorname{\rm H}}_{n-j-1}(Y_{a}^{b})$} ;
\draw (-5,1) node(051){$\widetilde{\operatorname{\rm H}}_{n-j-1}(Y_a)$} ;
\draw (-5,0) node(050){$ \operatorname{\rm H}^{j}(X_a)$} ;
\draw (0,-1) node(001){} (-3,-1) node(0301){$\operatorname{\rm H}^{j}(X_a^b)$}  ;

\draw  (1,2) node(12){$\widetilde{\operatorname{\rm H}}_{n-j-1}(Y_{b}^{c})$}  ;
\draw (-1,1) node(011){$\widetilde{\operatorname{\rm H}}_{n-j-1}(Y_b)$} (1,1) node(11){} (3,1) node(31){$\widetilde{\operatorname{\rm H}}_{n-j-1}(Y_c)$};
\draw (-1,0) node(010){$ \operatorname{\rm H}^{j}(X_b)$}  (3,0) node(30){$\operatorname{\rm H}^{j}(X_c)$};
\draw (0,-1) node(001){} (1,-1) node(101){$\operatorname{\rm H}^{j}(X_b^c)$} ;

\draw (5,2) node(52){$\widetilde{\operatorname{\rm H}}_{n-j-1}(Y_{c}^{d})$} ;
\draw (7,1) node(71){$\widetilde{\operatorname{\rm H}}_{n-j-1}(Y_d)$} ;
\draw (7,0) node(70){$ \operatorname{\rm H}^{j}(X_d)$} ;
\draw (5,-1) node(501){$\operatorname{\rm H}^{j}(X_c^d)$} ;

\draw[->] (051) to node[above]{} (032); 
\draw[->] (050) to node[right]{\rotatebox[origin=c]{90}{$\cong$}} (051);
\draw[->] (010) to node[right]{\rotatebox[origin=c]{90}{$\cong$}} (011);
\draw[->] (011) to node[above]{} (032); 
\draw[->] (0301) to node[above]{} (050);
\draw[->] (0301) to node[above]{} (010); 

\draw[->] (011) to node[above]{} (12); 
\draw[->] (101) to node[below]{} (010);
\draw[->] (101) to node[below]{} (30); 
\draw[->] (30) to node[right]{\rotatebox[origin=c]{90}{$\cong$}} (31);
\draw[->] (31) to node[below]{} (12); 

\draw[->] (31) to node[above]{} (52); 
\draw[->] (71) to node[below]{} (52);
\draw[->] (501) to node[below]{} (70); 
\draw[->] (70) to node[right]{\rotatebox[origin=c]{90}{$\cong$}} (71);
\draw[->] (501) to node[below]{} (30); 

\end{tikzpicture}
\end{center}
Applying the Diamond Principle to each of these diamonds, the four indecomposable summands change as follows:
$$
\begin{array}{rclcccccc}
\operatorname{\rm H}_j(\XX_{\{a, b, c, d\}}) & &  \widetilde{\operatorname{\rm H}}_{n-j-1}(\YY_{\{a, b, c, d\}}) \\
 & &   \\
\scalebox{0.065}{\input{dd.tex}} & \leftrightarrow & \scalebox{0.065}{\input{uu.tex}}\\
\scalebox{0.065}{\input{du.tex}} & \leftrightarrow &\scalebox{0.065}{\input{ud.tex}}\\
\scalebox{0.065}{\input{ud.tex}} &\leftrightarrow &\scalebox{0.065}{\input{du.tex}} \\
\scalebox{0.065}{\input{uu.tex}} &\leftrightarrow &\scalebox{0.065}{\input{dd.tex}}
\end{array}
$$
From here we conclude that
$$
\begin{array}{lcccccl}
{}_{n-j-1}\widetilde{\mu}_\YY^{{ \backslash\! \backslash}}(R) &=& \langle \scalebox{0.065}{\input{dd.tex}}\, |\, \widetilde{\operatorname{\rm H}}_{n-j-1}(\YY_{\{a, b, c, d\}}) \rangle & = & \langle \scalebox{0.065}{\input{uu.tex}}\, |\, \operatorname{\rm H}_j(\XX_{\{a, b, c, d\}}) \rangle
 & = & {}_j\mu_\XX^{^{/\!/}}(R),  \\
{}_{n-j-1}\widetilde{\mu}_\YY^{\vee}(R) &=&
\langle \scalebox{0.065}{\input{du.tex}}\, |\, \widetilde{\operatorname{\rm H}}_{n-j-1}(\YY_{\{a, b, c, d\}}) \rangle &=&
\langle \scalebox{0.065}{\input{ud.tex}}\, |\, \operatorname{\rm H}_j(\XX_{\{a, b, c, d\}}) \rangle
&=&  {}_j\mu_\XX^{^{\wedge}}(R) ,  \\
{}_{n-j-1}\widetilde{\mu}_\YY^{\wedge}(R) &= &
\langle \scalebox{0.065}{\input{ud.tex}}\, |\, \operatorname{\rm H}_{n-j-1}(\YY_{\{a, b, c, d\}}) \rangle  &=&
\langle \scalebox{0.065}{\input{du.tex}}\, |\, \operatorname{\rm H}_j(\XX_{\{a, b, c, d\}}) \rangle &=&
{}_j\mu_\XX^{^{\vee}}(R),\\
{}_{n-j-1}\widetilde{\mu}_\YY^{^{/\!/}}(R) &= &
\langle \scalebox{0.065}{\input{uu.tex}}\, |\, \widetilde{\operatorname{\rm H}}_{n-j-1}(\YY_{\{a, b, c, d\}}) \rangle&=&
\langle \scalebox{0.065}{\input{dd.tex}}\, |\, \operatorname{\rm H}_j(\XX_{\{a, b, c, d\}}) \rangle &=&
{}_j\mu_\XX^{^{\backslash\! \backslash}}(R).
\end{array}
$$
Since the measures are the same on all the rectangles, by the Equivalence Theorem the associated decorated diagrams are also the same. This proves Theorem \ref{theorem}.
\end{proof}

\begin{Example}\label{ex}
Revisiting Example \ref{example}, let $c_X$ denote the 0-dimensional homology cycle in $\XX$ that ceases to exist beyond $a_1$ and is killed at $a_2$ (see the picture below). 
\begin{figure}[h!]\label{cikel}
\begin{center}
\scalebox{0.5}{ \input{examplecycle.tex} } 
\end{center}
\end{figure}
The parametrized homology of $\XX$ and the reduced parametrized homology of $\YY$ can be seen in Figure~\ref{Diagramprimerjava}. The red point in the diagram representing the cycle $c_X$ indicates its dimension, whereas the symbol $\backslash\!\backslash$ designates the way it perishes at endpoints. By Theorem \ref{theorem} we know that there is a corresponding 1-dimensional homology cycle $c_Y$ in $\YY$. Its dimension is indicated in blue. In addition to the change of dimension, we observe the following: 
\begin{itemize}
\item
Since cycle $c_Y$ persists over $[a_1, a_2)$ like $c_X$, the decorations of the points representing these two cycles are the same;
\item
In contrast to $c_X$, the cycle $c_Y$ is killed at $a_1$ and ceases to exist beyond $a_2$. This is expressed by the symbol $/\!/$.
\end{itemize}
\begin{figure}[h!]
\begin{center}
\includegraphics[scale=0.23]{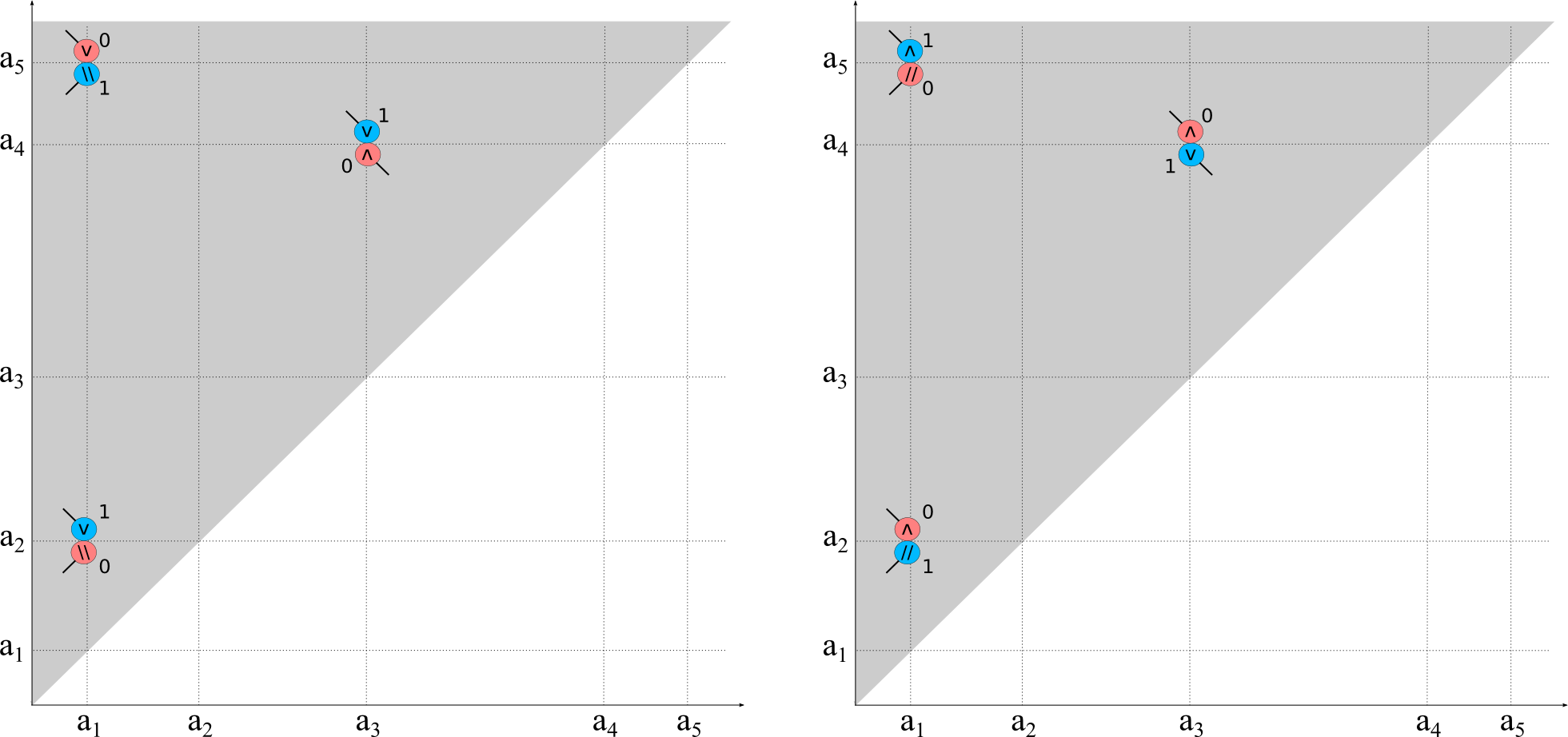}
 \end{center}
\caption{\label{Diagramprimerjava} The parametrized homology of $\XX$ is on the left and the one of $\YY$ on the right. }
\end{figure}
\end{Example}
\newpage
\paragraph{Acknowledgement}
This paper has benefited  greatly from the numerous discussions with Gunnar Carlsson, Vin de Silva, and Jaka Smrekar. The author also thanks Peter J. Verov\v{s}ek and Henry Adams for proofreading.

\printbibliography

\end{document}